\documentclass[sn-mathphys,Numbered]{sn-jnl}

 \usepackage{marginnote}

\usepackage{dsfont}
\usepackage{graphicx}%
\usepackage{multirow}%
\usepackage{amsmath,amssymb,amsfonts}%
\usepackage{amsthm}%
\usepackage{mathrsfs}%
\usepackage[title]{appendix}%
\usepackage{xcolor}%
\usepackage{textcomp}%
\usepackage{manyfoot}%
\usepackage{booktabs}%
\usepackage{algorithm}%
\usepackage{algorithmicx}%
\usepackage{algpseudocode}%
\usepackage{listings}%
\usepackage[center]{caption}
\usepackage{amsmath}
\newcommand*{\defeq}{\stackrel{\text{def}}{=}}
\def\eps{\epsilon}

\theoremstyle{thmstyleone}%
\theoremstyle{thmstyletwo}%

\theoremstyle{thmstylethree}%

\usepackage{titlesec} 

\titleformat{\subsubsection}[runin]
  {\normalfont\normalsize\bfseries}{\thesubsubsection}{1em}{}

\begin{document}

\title[A Three-dimensional tumor growth model and its boundary instability]{A Three-dimensional tumor growth model and its boundary instability}


\author[1]{\fnm{Jian-Guo} \sur{Liu}}\email{jian-guo.liu@duke.edu}
\equalcont{These authors contributed equally to this work.}

\author[2]{\fnm{Thomas} \sur{Witelski}}\email{witelski@math.duke.edu}
\equalcont{These authors contributed equally to this work.}

\author[3]{\fnm{Xiaoqian} \sur{Xu}}\email{xiaoqian.xu@dukekunshan.edu.cn}
\equalcont{These authors contributed equally to this work.}

\author*[4]{\fnm{Jiaqi} \sur{Zhang}}\email{jz374@duke.edu}

\affil[1,2]{\orgdiv{Mathematics Department}, \orgname{Duke University},  \state{North Carolina}, \country{US}}

\affil[3,4]{\orgdiv{Applied Mathematics and Computational Sciences}, \orgname{Duke Kunshan University},  \state{Jiangsu}, \country{China}}


\abstract{In this paper, we investigate the tumor instability by employing both analytical and numerical techniques to validate previous results and extend the analytical findings presented in a prior study by Feng et al 2023.
Building upon the insights derived from the analytical reconstruction of key results in the aforementioned work in one dimension (1D) and two dimensions (2D), we extend our analysis to three dimensions (3D). Specifically, we focus on the determination of boundary instability 
using perturbation and asymptotic analysis along with spherical harmonics.
Additionally, we have validated our analytical results in a two-dimensional framework by implementing the Alternating Directional Implicit (ADI) method, as detailed in Witelski and Bowen (2003).  Our primary focus has been on ensuring that the numerical simulation of the propagation speed aligns accurately with the analytical findings. Furthermore, we have matched the simulated boundary stability with the analytical predictions derived from the evolution function, which will be defined in subsequent sections of our paper. These alignment is essential for accurately determining the stability or instability of tumor boundaries. 
}




\maketitle

\section{Introduction}\label{sec1}

Cancer is one of the major concerning diseases around the world, and the mathematical study of tumor growth has proved to be meaningful in providing insights for treatments. Therefore, investigating partial differential equation (PDE) models in the context of tumor growth is of great importance \cite{cristini2003nonlinear,friedman2001symmetry,greenpan1972models}. 

Recent advancements in tumor growth modeling have emphasized nonlinear dynamics, boundary behaviors, and stability under perturbations. Cristini et al. (2003) identified key dimensionless parameters controlling tumor evolution in \cite{cristini2003nonlinear}, delineating growth models influenced by vascularization. Greenspan (1972) provided a foundational diffusion-based free boundary model for solid tumor growth in \cite{greenpan1972models}, while Friedman and Reitich (2000) revealed numerous asymmetric solutions in stationary tumor models in \cite{friedman2001symmetry}. These highlight the interaction between mathematical parameters and tumor morphology. Moreover, Feng et al. examined tumor boundary stability in \cite{feng2023tumor}, finding that tumor boundaries in vitro remain stable under various conditions, while in vivo boundaries exhibit instability when the nutrient consumption rate exceeds certain thresholds, influenced by the wave numbers of domain perturbations. Improved numerical scheme for tumor growth PDE models include the novel prediction-correction reformulation that accurately approximate the front propagation of speed with strong nonlinearity introduced by Liu et al. (2018) in \cite{liu2018accurate}. 
These studies underline the complex interplay between biological and mathematical factors in tumor growth, guiding our research on the stability of tumor growth models after perturbation in three dimensions.

The main question that we are focusing on is the instability of the tumor boundary induced by nutrient consumption and supply, which is significant since previous research has
suggested the shape of tumors as one of the crucial criteria to determine a tumor is benign or
malignant \cite{araujo2004history,byrne2006modelling,lowengrub2009nonlinear,roose2007mathematical}. 
Specifically, we investigate this based on \cite{feng2023tumor} and we use finite difference numerical schemes to match up to the two-dimensional (2D) analytical results presented in \cite{feng2023tumor}, as well as extending the results to three-dimensional perturbed spherical scenario. We first apply the Alternating Direction Implicit (ADI) numerical scheme as in Witelski (2003), to validate  against the two dimensional (2D) analytical results in Feng (2023) on tumor dynamics. Following this, we introduce spherical harmonics as a perturbation method, expanding our tumor model to a three-dimensional (3D) framework. In this extension, we employ modified Bessel functions, hyperbolic trigonometric functions, and asymptotic expansions to derive analytical results for the 3D scenario. 
This study is significant from different perspectives. Firstly, the findings reported in \cite{feng2023tumor} are predominantly based on analytical techniques, and the introduction of numerical simulations offers an alternative perspective on the 2D tumor PDE model, shedding light on boundary instability. Moreover, the extension of our investigation into 3D tumor dynamics, both analytically and numerically, has practical implications, providing deeper insights into the malignant potential of tumors.

In this paper, we begin by describing the model, detailing both in vitro and in vivo nutrient models. We then focus on 2D numerical schemes, matching our numerical simulations to the results in \cite{feng2023tumor}. Subsequently, we investigate boundary stability and introduce our analytical approach for constructing 3D scenarios, highlighting the distinction between our 2D and 3D tumor boundary instability. The paper concludes with a summary of our findings and a discussion of future research directions.

\section{Model formulation}\label{sec2}
\begin{figure}[!htb]
    \centering
    \includegraphics[scale=0.18]{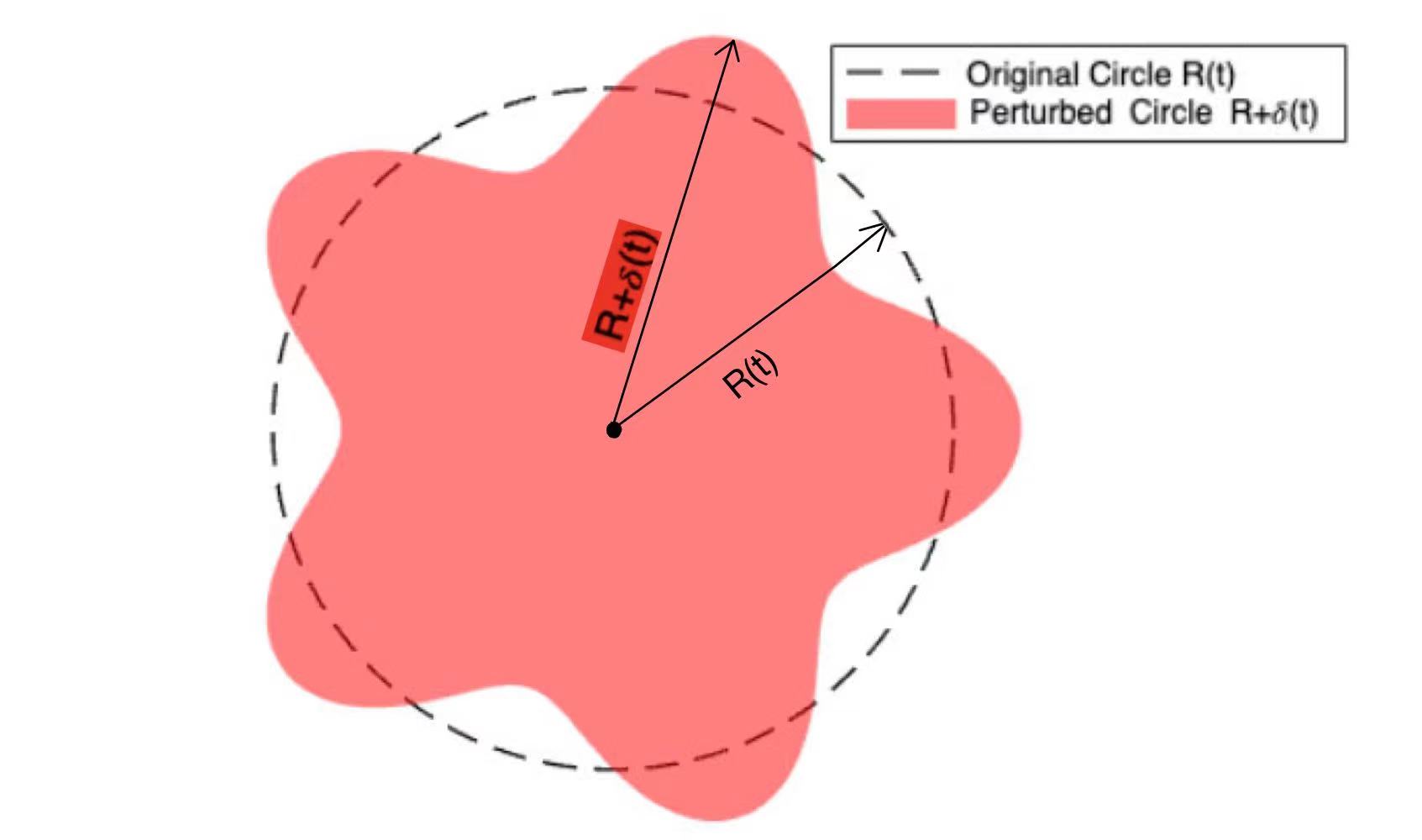}
     \captionsetup{justification=raggedright} 
     \captionsetup{justification=justified}
    \caption{Schematic figure of the tumor and the surrounding exterior. The cell population density \(\rho\) is defined in the red shaded tumoral region; \(c^{(o)}\) and \(c^{(i)}\) are the nutrient concentration defined outside and inside the red shaded tumoral region \cite{feng2023tumor}.  } 
    \label{schematic}
\end{figure}
The full PDE model in paper \cite{feng2023tumor} we are focusing on is based on the sharp interface limit of a cell density model. 
In Figure \ref{schematic}, we present a schematic representation of the tumor's geometry. This depiction illustrates a 2D circular domain with radius $R$, encapsulated within nutrient exterior. This simplified representation serves as a conceptual aid; the actual geometry of the tumor, governed by the equations detailed subsequently, can be any general domain. Our focus on a basic geometric form, such as the axisymmetric circle in 2D or the sphere in 3D (as explored in Section 5) is instrumental. It allows us to distill the essential physics of the model. By employing these simplified shapes, specifically the spherically symmetric sphere in our 3D analysis, we aim to elucidate the fundamental principles underlying the model's behavior.

 Following this, we are interested in physical properties such as the tumor density, the nutrient concentration and the pressure. Specifically, we denote $\rho(x,t)$ as the cell population density, and $c(x,t)$ as the nutrient concentration. The production rate of tumor cells $G_0$, is only dependent on the nutrient concentration, where we use $G(c)=G_0 c$ as growth function in this paper for simplicity. To represent the tumoral region at time $t$, we use $D(t)=\{\rho(x,t)>0\}$ to denote the support of $\rho$. Here, we assume the tumoral region expands with a finite speed governed by the Darcy's law $v=-\nabla p \vert_{\partial \Omega}$ via the pressure $p(\rho)={\rho}^k$ \cite{feng2023tumor}. In these models, we use Darcy's law since the flow within the tumor or biological tissue can be highly nonlinear, meaning that small changes in density can lead to large changes in flow or vice versa. A power law is a simple way to capture such nonlinear effects. The positive power $k>0$ ensures the consistency with physical laws, since the pressure in a fluid or gas typically increases as the density increases. Also, a positive power $k$ in the equation ensures that the pressure $p(\rho)$ is a well-defined, real-valued function for all positive values of density $\rho$. Furthermore, $k>0$ in the pressure $p(\rho)=\rho^k$ ensures the support of $\rho$ is compact. 
 Therefore, the density $\rho$ is governed by the Porous Medium Equation (PME) type equation:
\begin{equation} \label{govern}
    \frac{\partial \rho}{\partial t}  -\nabla \cdot (\rho \nabla p(\rho))=G_0 \rho c, \ \ \ x \in \mathbb{R}^n, \ \ \  t \geq 0,
\end{equation}
where $n=2$ or $n=3$. 

The consumption, exchange and diffusion of the nutrient are governed by the following reaction-diffusion equation for nutrient concentration $c$:
\begin{equation} \label{govern c}
    \tau {\partial c\over \partial t}- \nabla^2 c+ \Psi(\rho,c)=0
\end{equation}
where $\tau \ll 1$ is the characteristic time scale of nutrient evolution and $\Psi(u,c)$ describes the overall effects of the nutrient supply regime outside the tumor and the nutrient consumption by cells inside the tumor ~\cite{feng2023tumor}. 

Specifically, the overall function $\Psi(\rho,c)$ is mathematically written as $ \Psi(\rho,c)=\lambda \rho c \cdot \mathds{1}_{D}-\Tilde{\lambda}(c_B-c) \cdot \mathds{1}_{\partial D(t)} $. 
For simplicity, we set $\Tilde{\lambda} = 1$ and $\lambda >0$. $\mathds{1}_{D}$ indicates the indicator function of set D (the region of the tumor), and $\mathds{1}_{\partial D(t)}$ is the indicator function of the complement of set $D$ \cite{feng2023tumor}, where $\partial D(t)$ relates to the evolving tumor boundary. The indicator function of a set $A$ is defined such that it assigns a value of 1 to an element $x$ if $x$ is in the set $A$, and a value of 0 if $x$ is not in the set $A$.

For simplicity in later analysis, we consider the following elliptic formulation by dropping the time derivative, which is reasonable since $\tau \ll 1$ :
\begin{equation} 
    -\nabla^2 c+ \Psi(\rho,c)=0
    \label{govern real}
\end{equation}
This further yields:
\begin{equation} \label{sect2 key c}
    \begin{aligned}
        -\nabla^2 c+\lambda \rho c = 0, \ \ \mbox{for} \ \ x \in D(t)
    \end{aligned}
\end{equation}
where $n=2,3$, $D(t)$ is the tumor region, and $\lambda$ is the consumption rate \cite{feng2023tumor}.

Taking the incompressible limit of $k \to \infty$, and consider the patch solutions $\rho_\infty=\mathds{1}_{D_\infty}$, $D_\infty \subset \mathbb{R}^n$, yield the simplified versions of equations (\ref{govern}) and (\ref{govern real}):
\begin{eqnarray} \label{p real}
    -\nabla^2 p_{\infty}=G_0 c \ \ \mbox{for} \ \ x \in D_{\infty}(t), \\ \nonumber
    p_{\infty}=0, \ \ \mbox{for} \ \ x \in \mathbb{R}^n \backslash D_{\infty}(t).
\end{eqnarray}
 The boundary conditions for pressure $p$ in the in vitro model is that pressure is bounded everywhere and the pressure value and its derivative are zero at the boundary. Therefore, these give us the complete PDE problem for solving the pressure $p$ for the in vitro model.
\begin{equation} \label{c real}
     -\nabla^2 c+\lambda c = 0, \ \ \mbox{for} \ \ x \in D_\infty(t)
\end{equation} 
where $n=2,3$, and the in vitro and in vivo models will be introduced shortly in section 2.2 \cite{feng2023tumor}. Combining (\ref{p real}) and (\ref{c real}), we can denote
\begin{equation} \label{useful combo}
u= p+\frac{G_0}{\lambda} c\qquad \mbox{satisfying}\quad
    \nabla^2\left(p+\frac{G_0}{\lambda}c\right)=0\,.
\end{equation}
This gives us the Laplace's equation for $u$, which will be useful for later constructions.

The domain we are considering here is the disk or sphere with radius $R$, where the solution is bounded everywhere. The boundary conditions and initial conditions for the above PDEs will be introduced for the simplified versions of this full models in the following subsections. Further, the boundary moving speed in the direction normal to the boundary at a point $x$, denoted by $\sigma(x)$, is given by:
\begin{equation}
    \sigma(x)=-\nabla p_{\infty} \cdot \hat{n}(x)
\end{equation}
where $\hat{n}(x)$ is the outer unit normal vector at $x\in \partial D_{\infty}(t)$. 

\subsection{Two models for the nutrient distribution}
Following the general PDEs model in Equation~(\ref{c real}) and Equation~(\ref{govern real}), we consider the two simplified versions with different concentration assumptions, as introduced in \cite{feng2023tumor}, which are the in vitro and in vivo models. 
\subsubsection{The in~vitro model.} 
For the in vitro version, it is assumed that the tumor is surrounded by a liquid and the exchange rate with the background is extremely high so that $c$ matched $c_B$ on the tumor boundary \cite{feng2023tumor}. This gives us constant nutrient concentration outside the tumor and quasisteady inside the tumor. 
The in vitro concentration assumption outside of tumoral region coupled with (\ref{c real}) give us the system of equations for the in vitro regime \cite{feng2023tumor}:
\begin{equation}  \label{in vitro vitro}
    \begin{aligned}
        -\nabla^2 c+\lambda c = 0, \ \ \mbox{for} \ \ x \in D_\infty(t), \\
         c=c_B, \ \ \mbox{for} \ \ x \in \mathbb{R}^n \backslash D_\infty(t),
    \end{aligned}
\end{equation}
where $n=2,3$ and $\lambda>0$ is the consumption rate.
Therefore, (\ref{in vitro vitro}) together with the boundary condition $c(R)=c_B$ gives us the complete PDE problem for solving the concentration $c$ for the in vitro regime model. We only need to solve the concentration $c$ inside the tumor since $c$ is always constant outside the tumor. 

\subsubsection{The in~vivo model.} 
For the in vivo version, the nutrient is transported by vessels outside the tumor and reaches $c_B$ at the far field. Correspondingly, it is assumed that the exchange rate outside the tumor is determined by the concentration difference from the background, i.e., $c_B-c$ \cite{feng2023tumor}. The in vivo concentration assumption outside of tumoral region coupled with (\ref{c real}) give us the system of equations for the in vivo regime \cite{feng2023tumor}:
\begin{equation} \label{govern vivo}
   \begin{aligned}
        -\nabla^2 c+\lambda c = 0, \ \ \mbox{for} \ \ x \in D_\infty(t), \\
        -\nabla^2 c=c_B-c, \ \ \mbox{for} \ \ x \in \mathbb{R}^n \backslash D_\infty(t),
  \end{aligned}    
\end{equation}
where $n=2,3$, and the concentration is quasisteady inside and outside the tumor. 
If we denote the concentration inside the tumor as $c^{(i)}(r,t)$, and the concentration outside the tumor as $c^{(o)}(r,t)$, then the solution for concentration is everywhere bounded, and boundary conditions for concentration $c$ in  the in vitro model is
\begin{equation} \label{c vivo bcs}
    \begin{aligned}
        c^{(i)}(R,t)= c^{(o)}(R,t) \\     
        \frac{\partial c^{(i)}(r,t)}{\partial r} \bigg|_{r=R} =  \frac{\partial c^{(o)}(r,t)}{\partial r} \bigg|_{r=R}
    \end{aligned}
\end{equation}
Therefore, (\ref{govern vivo}) and (\ref{c vivo bcs}) give us the complete PDE problem for solving the concentration $c$ for the in vivo regime model. 

The PDE problem for pressure \( p \) in the in vivo model is identical to that in the in vitro model. Therefore, equation (\ref{p real}) along with the boundary conditions outlined prior to Section 2.1, comprehensively define the PDE problem for determining pressure \( p \) in both the in vitro and in vivo models.

\section{The Two-Dimensional 
Problem}
In this section, we will examine the 2D analytical findings on tumor boundary stability and instability as presented in \cite{feng2023tumor}, where the focus is on radially symmetric boundary. In 2D, the tumoral region is considered as a disk with radius $R$ for simplicity. This aims to establish the foundational equations pertinent to our research and facilitate a comparative analysis with the 3D results we have obtained. Specifically, in section 3.1, we will skip their calculation details and provide a comprehensive overview of their procedural steps. Subsequently, in Section 3.2, we will state their results for 2D.
\subsection{Overview of the steps}
In both the in vitro and in vivo models, similar methodologies are employed, yet they yield distinct results. In this section, we will primarily focus on the more complex in vivo scenario as an illustrative example. The results for both the in vitro and in vivo models will be presented in Section 3.2. \\ \\
\textbf{Step 1. Solve radial symmetric solution.} \\
    In \cite{feng2023tumor}, the authors solve the nutrient concentration using the second kind modified Bessel functions $I_n(r)$ and $K_n(r)$. Then they solve the pressure $p_\infty$ and use the symmetry to derive $p_\infty(r,t)$.
The velocity of the boundary for the in vitro model is 
\begin{equation}
    \frac{dR}{dt}=-\partial_r p(R(t))=\frac{G_0 c_B I_1(\sqrt{\lambda}R)}{\sqrt{\lambda} I_0(\sqrt{\lambda}R)}
\end{equation}
and for the in vivo model:
    \begin{equation}
        \frac{dR}{dt}=-\partial_r p(R(t))
        =\frac{G_0 c_B K_1(R) I_1(\sqrt{\lambda}R)}{\lambda K_0(R) I_1(\sqrt{\lambda}R)+\sqrt{\lambda}K_1(R)I_0(\sqrt{\lambda}R)}
    \end{equation}
\textbf{Step 2. Define the perturbation in terms of parameter curve.} \\
At time $t$, the authors introduce a perturbation characterized by an amplitude $\delta(t)$ and the boundary of this perturbation is described as follows:
\begin{equation} \label{perturbation}
    \Tilde{B_t}(\theta)=\{(r,\theta)| r=R(t)+\delta(t)\cos(m \phi), \phi \in [-\pi, \pi)\}.
\end{equation}
\textbf{Step 3. Expression of asymptotic solutions on the perturbed region.}\\
The perturbed solutions (drop the tilde) $c$ and $p$ have the asymptotic expansions:
\begin{equation}
    \begin{aligned}
        c(r, \phi, t)=c_0(r,t)+\delta(t)c_1(r, \phi, t)+O(\delta^2) \\
         p(r, \phi, t)=p_0(r,t)+\delta(t)p_1(r, \phi, t)+O(\delta^2) 
    \end{aligned}
\end{equation}
where the leading order terms $c_0$ and $p_0$ correspond to the unperturbed solutions solved in Step 1. The main response corresponding to the perturbation are captured by the first-order terms $c_1(r, \theta, t)$ and $p_1(r, \theta, t)$. Since the boundary conditions imply that  all terms subsequent to the leading term in the series expansion for the concentration $c$ are identically zero, we can simplify the expansions and approximate them as:
\begin{equation}
    \begin{aligned}
      c_1(r, \phi, t) \approx c_1(r,t) \cos(m \phi), \\ 
      p_1(r, \phi, t) \approx  p_1(r,t) \cos(m \phi).
    \end{aligned}
\end{equation}
\textbf{Step 4. Equations for the first order terms}\\
Plug the asymptotic solutions into the model on $\Tilde{D}(t)$ and solve for $c_1^{(i)}$ and $c_1^{(o)}$:
\begin{equation}
    \begin{aligned}
        -\nabla^2 c_1^{(i)} (r, \phi,t)+\lambda c_1^{(i)} (r, \phi,t)=0, \\ 
        -\nabla^2 c_1^{(o)} (r, \phi,t)+ c_1^{(o)} (r, \phi,t)=0.
    \end{aligned}
\end{equation}
On the other hand, $p_1$ solves:
\begin{equation}
    \nabla^2 p_1(r, \phi, t)=G_0 c_1^{(i)}(r,\phi, t), \ \ \mbox{in} \ \ \Tilde{D} (t).
\end{equation}
\textbf{Step 5. Determine the coefficients by matching the boundary}\\
The authors evaluate the perturbed solution on the boundary point $\Tilde{B}(R+\delta \cos(m \phi), \phi)$ by first taking Taylor expansion around $R$ with respect to the first variable up to first order:
\begin{equation} \label{F4}
    \begin{aligned}
        c^{(i)}(R+\delta \cos(m \phi), \phi)&=c^{(o)}(R+\delta \cos(m \phi), \theta) \\ 
        & =c_0(R)+\delta \partial_r c_0(R) \cos(m \phi)+\delta  c_1(R) \cos(m\phi) , \  \forall \theta \in [-\pi,\pi], \\ \\
        \partial_r c^{(i)}(R\!+\!\delta \cos(m \phi),\phi)&= \partial_r c^{(o)}(R+\delta \cos(m \phi), \phi)\\
    &\hspace{-1em}=\partial_rc_0(R)\!+\! \partial_r^2 c_0(R)\delta \cos(m \phi)\!+\! \delta  \partial_r c_1(R) \cos(m \phi),        
        \ \forall \phi \in [-\pi,\pi].
    \end{aligned}
\end{equation} 
Then by using the boundary condition for $c_1$, we require 
the solution to be bounded everywhere and $c(\infty, \theta)=c_B$ for any $\phi \in [-\pi, \pi)$ yields:
\begin{equation}
\begin{aligned}
    c_1^{(o)}(\infty,t)=0,
\end{aligned}
\end{equation}
Using the above boundary conditions, the coefficients $a_1(R(t))$ and $b_1(R(t))$ are determined.\\ \\
\textbf{Step 6. Determine the amplitude evolution equation $\delta^{-1} \frac{d \delta}{dt}$} for $p_1$. \\
The perturbed pressure solution $p$ satisfies the boundary condition $p=0$ at  $\Tilde{B}_t$. Evaluating $p$ at $\Tilde{B}(R+\delta\cos(m\phi)) \in \Tilde{B}_t$ yields
\begin{equation}
    p(R+\delta \cos(m \phi), \theta)=p_0(R)+\partial_r p_0(R) \delta \cos(m \phi) + p_1(R)\delta \cos(m \phi)+O(\delta^2).
\end{equation}
where  $O(1)$ yields $dR/dt= -\partial_r p_0(R)$, and $O(\delta)$ yields the linear $d\delta/dt$ equation. Since $p(R+\delta \cos(m \phi), \phi)=p_0(R)=0$, then
\begin{equation}
    \partial_r p_0(R)+p_1(R)=0.
\end{equation}
By using the expression of $c_1$ and $p_0$, the particular solution of $p_1$ can be determined. After this, plugging the expression of $p$ into the normal boundary speed (Darcy's law) and taking the Taylor expansion for the $r$ variable yields
\begin{equation}
    \frac{dR}{dt}+\frac{d\delta}{dt} \cos(m \phi)=-(\partial_r p_0(R)+\partial_r^2 p_0(R) \delta \cos(m \phi)+\partial_r p_1(R) \delta \cos(m \phi))+O(\delta^2).
\end{equation}
which simplifies to 
\begin{equation}
    \delta^{-1} \frac{d\delta}{dt}=-(\partial_r^2p_0(R)+\partial_r p_1(R))+O(\delta).
\end{equation}
This is used to determine the evolution of the perturbation magnitude by the sign of $\delta^{-1} \frac{d\delta}{dt}$. Specifically, if it is positive, then it implies the growing of the magnitude, which leads to instability with the growing perturbation. If it is negative, then it implies the decay of the perturbation amplitude and leads to stability of the tumor boundary.
\subsection{Results of 2D tumor boundary stability/instability in \cite{feng2023tumor}}
In conclusion, given $G_0>0, c_B>0, \lambda>0$, and $m \in \mathbb{N}^+.$ When the radius of the tumor is around $R$. For the \textbf{in vitro} model, $\delta^{-1} \frac{d\delta}{dt}$ is given by:
\begin{equation}
    \delta^{-1} \frac{d\delta}{dt}=\frac{G_0 c_B I_1(\sqrt{\lambda}R)}{I_0(\sqrt{\lambda}R)} \left(\frac{I_1^{'}(\sqrt{\lambda}R)}{I_1(\sqrt{\lambda}R)}-\frac{I_m^{'}(\sqrt{\lambda}R)}{I_m(\sqrt{\lambda}R)}\right) + O(\delta)=F_3(\lambda,m,R)+O(\delta)
\end{equation}
For the \textbf{in vivo} model, $\delta^{-1} \frac{d\delta}{dt}$ is given by:
\begin{equation}
\begin{aligned}
    \delta^{-1} \frac{d\delta}{dt}&=\frac{G_0 c_B m}{\sqrt{\lambda}RC(R)} \left(\frac{C_1(R)}{C_m(R)}K_m(R) I_m(\sqrt{\lambda}R)-K_1(R) I_1( \sqrt{\lambda}R) \right)   \\   
    & -\frac{G_0 c_B}{C(R)} \left( \frac{C_1(R)}{C_m(R)}K_m(R) I_m^{'}(\sqrt{\lambda}R)-K_1(R) I_1^{'}( \sqrt{\lambda}R) \right)+O(\delta)\\
    &=F_4(\lambda,m,R)+O(\delta)
\end{aligned}
\end{equation}
where $C(R)=\sqrt{\lambda}K_0(R)I_1(\sqrt{\lambda}R)+K_1(R)I_0(\sqrt{\lambda}R)$. $C_j(R)=K_{j}'(R)I_{j}(\sqrt{\lambda}R)-\sqrt{\lambda}I_{j}'(\sqrt{\lambda}R)K_j(R)$. Here $m$ denotes the wave number, $I_n(x)$ and $K_n(x)$ are the $n$th modified Bessel functions. 

\begin{figure}[htp] 
        \centerline{\includegraphics[width=1.05\linewidth]{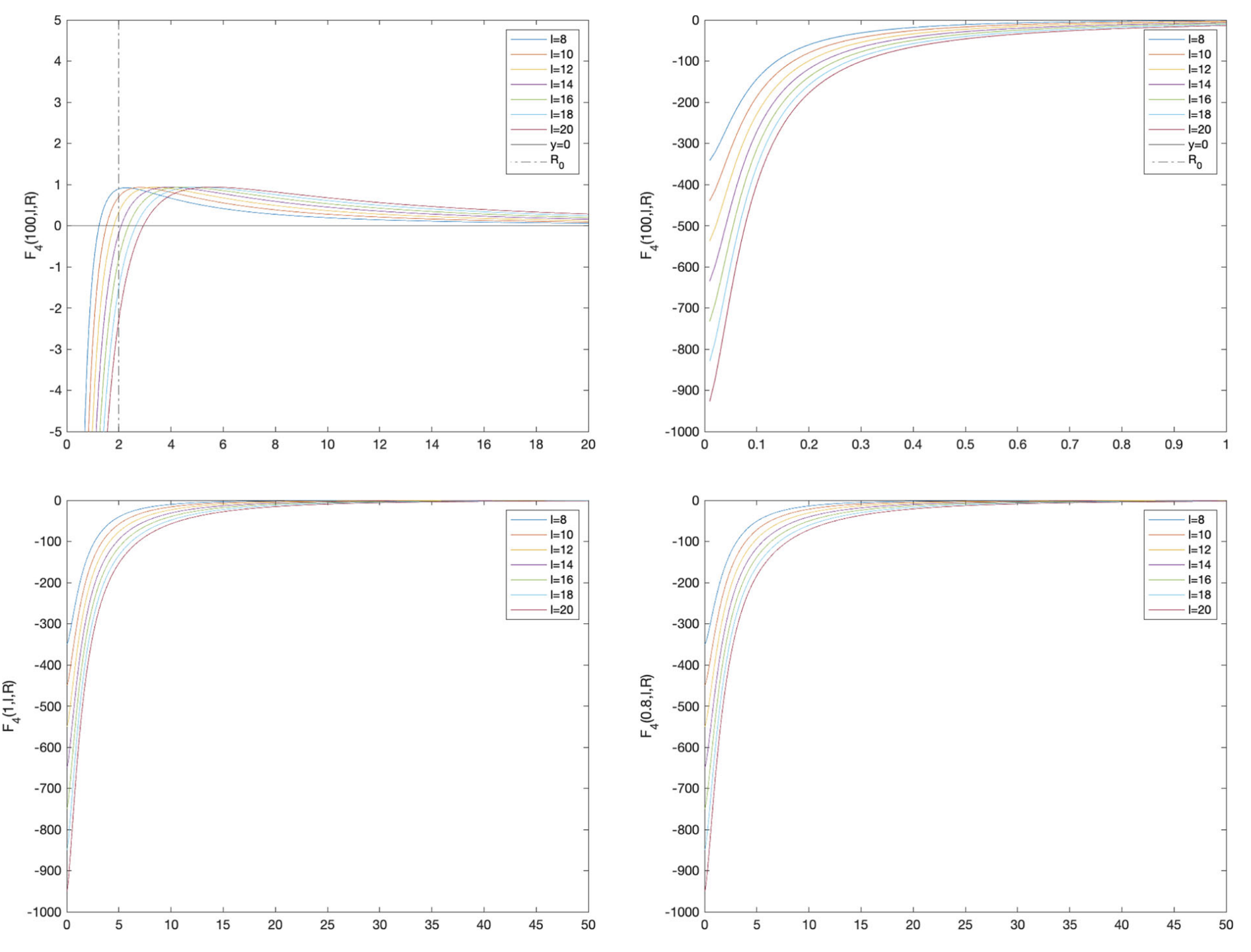}}
        \captionsetup{justification=justified}
        \caption{2D in vivo evolution function with $G_0=1, c_B=100$; top (left): $\lambda=100$ and $R \in [0,20]$; top (right): $\lambda=100$ and $R \in [0,1]$; bottom (left): $\lambda=1$ and $R\in[0.50]$; bottom (right): $\lambda=0.8$ and $R\in[0,50]$. Here the $l$ denotes as wave number. Figure from \cite{feng2023tumor}.} \label{2d results}
   \end{figure}

Therefore, the sign of the leading order terms of $\delta^{-1} \frac{d\delta}{dt}$ (such as $F_3(\lambda,m,R)$ and $F_4(\lambda,m,R) $ as shown above) determines the sign of $\delta^{-1} \frac{d\delta}{dt}$, thus implies the stability/instability of the tumor boundary after perturbation. The plots of the evolution function for the in vivo regime is shown in Figure~\ref{2d results}. The simulations of the tumor boundary evolution using ADI method are shown in Figure~\ref{3}, \ref{4}, where the detailed descriptions for the numerical scheme will be presented in Appendix~\ref{secA1}.
\begin{figure}[htp]
    \includegraphics[scale=0.38]{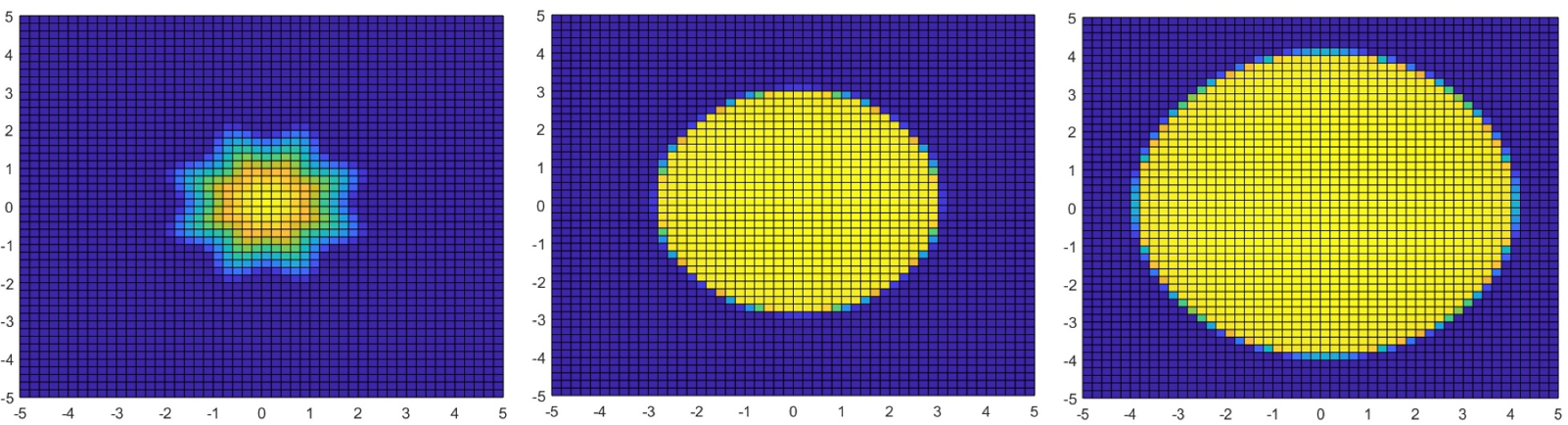}
    \captionsetup{justification=justified}
    \caption{ 2D tumor growth boundary evolution with $\lambda=0.5$ (three time profiles evolving from left to right) for an initial perturbation with wave number $m=8$} \label{3}
    \includegraphics[scale=0.38]{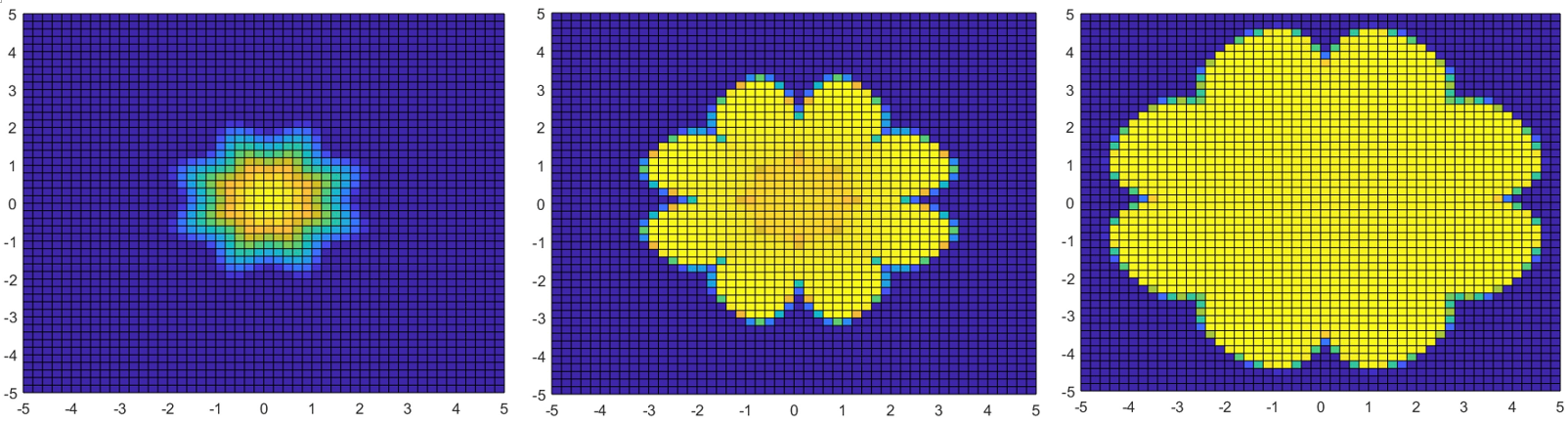}
    \captionsetup{justification=justified}
    \caption{ 2D tumor growth boundary evolution with $\lambda=40$ (three time profiles evolving from left to right) for an initial perturbation with wave number $m=8$}   \label{4}
\end{figure} \\ \\
\textbf{Remark 3.1.}
In the in vitro regime, both graphical and analytical analysis indicate that the growth of the perturbation amplitude and boundary stability are maintained. This holds true regardless of variables such as the consumption rate \( \lambda \), perturbation wave number $m \in \mathbb{N}$, and tumor size \( R \).\\ \\
\textbf{Remark 3.2.} For the in vivo regime,\\
1. When $0<\lambda \leq 1$ ($\lambda$: nutrient consumption rate), the boundaries show the same qualitative behavior as the in vitro model.\\
2. When $\lambda>1$, for any fixed wave number $m \geq 2$, there exists a threshold $R^{*}(m)$ such that $F_4(\lambda,m,R) < 0$ for $0<R<R^{*}(m)$, and $F_4(\lambda,m,R)>0$ for $R>R^{*}(m)$. This is proved by checking the asymptotes of $F_4$. Also, fix a proper value of $R_0$, there exists $m_0$ such that $F_4(\lambda,m,R_0)>0$ for $m<m_0$ and $F_4(\lambda,m,R_0)<0$ for $m>m_0$. Lastly, as the tumor size expands, $R(t)$ exceeds more thresholds $R^{*}(m)$, therefore the corresponding wave number perturbation become unstable successively. 

\section{The Three-Dimensional 
Problem}
While the foundational configurations and methodologies in the 3D analysis are similar to those employed in the 2D case, the extension to a 3D framework in exploring the tumor growth model is of substantive significance. Primarily, this 3D approach yields results with higher applicability to real-world scenarios, enhancing the model's practical relevance. Additionally, the 3D model employs more intricate techniques, such as spherical harmonics, to capture the complexities of tumor growth. These advanced methodologies facilitate a deeper understanding of tumor boundary instability, a pivotal aspect of our study. The technical intricacies and the evolution function construction that elucidate the tumor boundary's instability are detailed in Sections 4.3 and 4.4.  

In this section, we focus on analyzing the tumor growth model in a 3D scenario, examining it under both in vitro and in vivo nutrient conditions. Specifically, we explore the spherically symmetric model in a 3D context for these nutrient models before any perturbation, as detailed in sections 4.1 (in vitro) and 4.2 (in vivo). We then shift our attention to the stability of the 3D spherically symmetric model under these models, discussed in sections 4.3 (in vitro) and 4.4 (in vivo). A key part of our investigation involves reconstructing critical variables—concentration (\(c\)), pressure (\(p\)), and boundary speed (\(\sigma\))—and analyzing their behaviors both before and after perturbations. 
We also derive the evolution function \(\delta^{-1}\frac{d \delta}{dt}\) and use it to assess the stability of the tumor boundary, based on its sign. The critical threshold for the transition between stability and instability is identified through asymptotic analysis of spherical modified Bessel functions and hyperbolic functions, which play a significant role in the evolution function. Detailed results and methodologies are included in each respective section.
In Section 4.5, we emphasize a comparative analysis between the in vitro and in vivo results, highlighting biological insights gained from this comparison.
\subsection{\textbf{3D spherically symmetric solution for the \textbf{in vitro} model}}\mbox{} \\
\leavevmode
For the 3D scenario, based on the models given in section 2, we can reduce the problem from Bessel functions to hyperbolic trigonometric functions $\cosh$ and $\sinh$, since it is in odd dimensions. The Laplacian in 3D spherical coordinates (spherically symmetric) becomes
\begin{equation} \label{3d}
    \nabla^2 u = \frac{1}{r^2} \frac{\partial}{\partial r} \left( r^2 \frac{\partial u}{\partial r}  \right)
\end{equation}
\textbf{Step 1: Construct the concentration $c$} \\ 
The corresponding governing equation for concentration $c$ in three-dimensional (3D) spherical coordinate (spherically symmetric) is Eq. (2.26) in \cite{feng2023tumor}:
\begin{equation} 
    \begin{aligned} 
        - \frac{1}{r^2} \partial_r(r^2 \partial_r c) + \lambda c = 0,\ \ \  \mbox{for} \ \ r \leq R(t),  \\
        c=c_B, \ \ \ \mbox{for} \ \ r \geq R(t).
    \end{aligned}
\end{equation}
Since $u(r, \theta, \phi)$ is independent on $\theta$ and $\phi$.
Expanding the right hand side of Equation~(\ref{3d}) yields
\begin{equation}
    \begin{aligned}
        \nabla^2 u =\frac{2}{r}\frac{\partial u}{\partial r} + \frac{\partial^2 u}{\partial r^2}      
    \end{aligned}
\end{equation}
This is shown to be the same as $\frac{1}{r} \left( \frac{\partial^2}{\partial r^2} (ru) \right) $ by expanding this out using product rule:
\begin{equation} \label{3d-2}
    \frac{1}{r} \left( \frac{\partial^2}{\partial r^2} (ru) \right)
    =\frac{2}{r}\frac{\partial u}{\partial r} + \frac{\partial^2 u}{\partial r^2}  
\end{equation}
Hence, we have 
\begin{equation}
    \frac{1}{r} \left( \frac{\partial^2}{\partial r^2} (ru) \right) = -\lambda u
\end{equation}
we let $v=ru$, then we obtain the general solution 
\begin{equation}
    \begin{aligned}
        \frac{\partial^2}{\partial r^2} (v) = -\lambda v \\
        \Rightarrow v = C_1 \sin(\sqrt{\lambda}r) + C_2 \cos(\sqrt{\lambda}r) \\ 
        \Rightarrow u = \frac{1}{r}( C_1 \sin(\sqrt{\lambda}r) + C_2 \cos(\sqrt{\lambda}r) )
    \end{aligned}
\end{equation}
where $C_1$ and $C_2$ are arbitrary constants, and this is applicable to odd dimensions. Further applying the boundary condition, we can solve this equation:
\begin{equation}
    c(r,t) = c_B \frac{R}{r} \frac{ \sinh(\sqrt{\lambda} r)}{\sinh(\sqrt{\lambda}R)}
    \label{c0vitro}
\end{equation}
\textbf{Step 2: Construct the pressure $p$.} \\ 
We have the PDE for pressure $p$:
\begin{equation}  \label{2.27pressure}
    \begin{aligned}
        - \frac{1}{r^2} \partial_r(r^2\partial_r p) = G_0 c,\ \ \  \mbox{for}\ \ r \leq R,  \\
        p=0, \ \ \ \mbox{for} \ \ r \geq R \\
        \frac{dR}{dt}=-\partial_r p(R) \\
        \partial_r p(0)=0.
    \end{aligned}
\end{equation}
and we can use $c(r,t)$ from \eqref{c0vitro} for $r \leq R$ obtained from step 1. 

We can combine the equations of concentration ($c$) for $r \leq R(t) $ and the equation of pressure ($p$) for $r \leq R(t) $, done similarly to section 2, to obtain the Laplace's equation for $u=p+\frac{G_0}{\lambda}c$ as shown in (\ref{useful combo}), we have:
\begin{equation} \label{combo-1}
    \frac{d^2u}{dr^2}+\frac{2}{r}\frac{du}{dr}=0
\end{equation}
This is a Cauchy-Euler (CE) equation, and thus we plug in the trial solution $u=r^t$, which gives 
\begin{equation}\label{now 67 C-E}
   u = A_1(R) + A_2(R) r^{-1}
\end{equation}
where $A_1(R)$ and $A_2(R)$ are arbitrary constants. Since we require $c, p$ to be bounded at the origin, $A_2(R)$ must be zero, so $u$ must be constant in space, $u=A_1(R)$. 
Or we can expand it by product rule, since it is in odd dimensions, and then we have the following derivation, which leads us to the same result:
\begin{equation}
    \begin{aligned}
        \frac{1}{r^2}(2r p_r + r^2 p_{rr})+\frac{G_0}{\lambda}(\frac{1}{r^2}(2r c_r + r^2 c_{rr})=0 \\
        \frac{1}{r^2}[2r(p_r+(\frac{G_0}{\lambda}c)_r) + r^2(p_{rr}+(\frac{G_0}{\lambda}c)_{rr})]=-\frac{1}{r} \frac{d^2}{dr^2}(rv)=0 \ \ \mbox{for} \ \ v= p+\frac{G_0}{\lambda}c \\
        \frac{d^2}{dr^2}(rv)=0 \\
        rv=A_1(R)r+A_2(R) \\
        v = A_1(R)+\frac{A_2(R)}{r} \\
        v = A_1(R) \  \mbox{to ensure no blow up}
    \end{aligned}
\end{equation}
Then we can apply the boundary condition from \eqref{2.27pressure}$_2$ and determine the constant to be

\begin{equation}
        p(r=R)=A_1(R)-\frac{G_0}{\lambda}c_B=0 \qquad\implies \qquad
        A_1(R) = \frac{G_0}{\lambda}c_B
\end{equation}
Therefore, we derive the final solution for pressure $p$ for the in vitro model:
\begin{equation}
    \label{p-3d-vitro}
        p(r,t) =\frac{G_0 c_B}{\lambda} \left( 1- \frac{R \sinh(\sqrt{\lambda} r)}{r\sinh(\sqrt{\lambda} R)}  \right) \,.
\end{equation}
\textbf{Step 3: Construct the boundary velocity.} \\
Substitute the solution for $p$ we derived in (\ref{p-3d-vitro}) into 
\begin{equation}
    \frac{dR}{dt}=-\frac{\partial p}{\partial r} (R(t)) 
\end{equation}
We get the boundary velocity for the 3D spherically symmetric in vitro model is 
\begin{equation} \label{3d vitro speed}
     \frac{dR}{dt} = G_0 c_B  \left(  - \frac{ 1}{ \lambda R}  + \frac{  \cosh(\sqrt{\lambda}R)}{ \sqrt{\lambda} \sinh(\sqrt{\lambda}R)}  \right) .
\end{equation}

\subsection{\textbf{3D spherically symmetric solution for the \textbf{in vivo} model}}\mbox{} \\
\textbf{Step 1: Construct the concentration $c$} \\ 
The equations for concentration in the in vivo model can be described as Eq.(2.31) in \cite{feng2023tumor}
\begin{equation}
    \begin{aligned}
        -\frac{1}{r^2} \partial_r (r^2\partial_r c)+\lambda c=0, \ \ \ \mbox{for} \ \ r \leq R(t) \\
        -\frac{1}{r^2} \partial_r (r^2\partial_r c) = c_B - c, \ \ \ \mbox{for} \ \ r \geq R(t)
    \end{aligned}
\end{equation}
We solve the above two equations and then apply the boundary conditions, which are from the continuity of both $c$ and $\partial_r c$ at $R(t)$
\begin{equation}
    \begin{aligned}
        c^{(i)}(R,t) = c^{(o)}(R,t) \\ 
        \partial_r c^{(i)}(R,t) = \partial_r c^{(o)}(R,t)
    \end{aligned}
\end{equation}
and we get
\begin{equation} \label{also 3d c}
    c(r,t) = 
    \begin{cases}
        c_B a_0(R) \sinh(\sqrt{\lambda} r)/r \defeq c^{(i)}(r,t) &\mbox{for} \ \ r \leq R(t), \\
        c_B(b_0(R)e^{-r}+r)/r \defeq c^{(o)}(r,t) & \mbox{for} \ \ r \geq R(t),
    \end{cases}
\end{equation}
where $a_0$ and $b_0$ are given by:
\begin{eqnarray} \label{c-3d-vivo}  
       a_0(R) = \frac{R+1}{\sqrt{\lambda}\cosh(\sqrt{\lambda}R)+\sinh(\sqrt{\lambda}R)} \\ 
       b_0(R) = -\frac{\sqrt{\lambda}\, R \cosh(\sqrt{\lambda}R)-\sinh(\sqrt{\lambda}R)}{e^{-R}(\sqrt{\lambda}\cosh(\sqrt{\lambda}R)+\sinh(\sqrt{\lambda}R))}
\end{eqnarray}
\textbf{Step 2: Construct the pressure $p$} \\
We have the PDE for pressure $p$:
\begin{equation}  \label{2.27pressure}
    \begin{aligned}
        - \frac{1}{r^2} \partial_r(r^2\partial_r p) = G_0 c,\ \ \  \mbox{for}\ \ r \leq R(t),  \\
        p=0, \ \ \ \mbox{for} \ \ r \geq R(t) \\
        \frac{dR}{dt}=-\partial_r p(R) \\ 
        \partial_r p(0)=0.
    \end{aligned}
\end{equation}
as in (\ref{2.27pressure}), and we have the solution we get for concentration $c$ as in equations (\ref{also 3d c}) and (\ref{c-3d-vivo}).

\begin{figure}[htp] 
    \begin{minipage}{0.53\textwidth}
        \centering
        \includegraphics[width=1.0\linewidth]{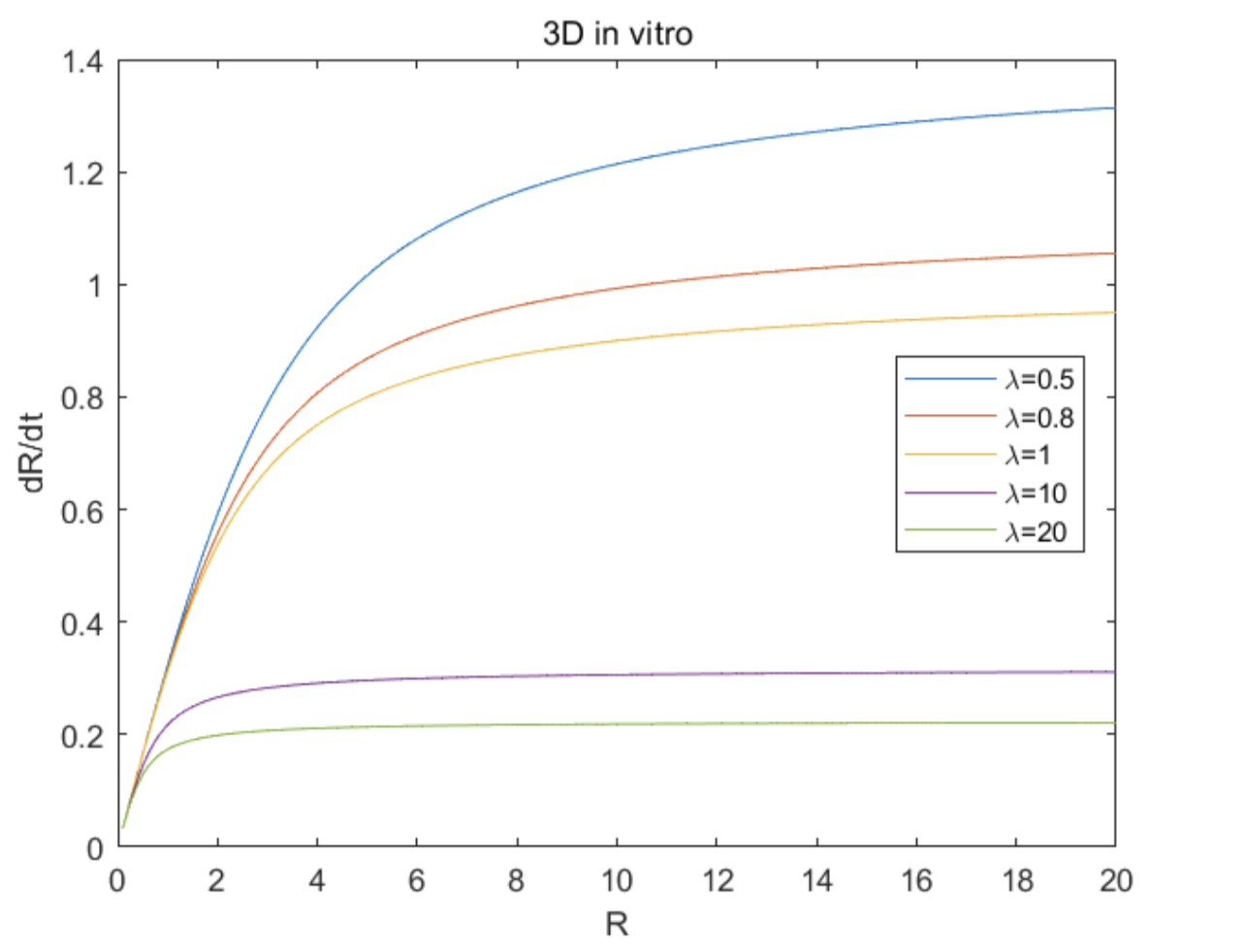}
        \caption{3D in vitro tumor evolution (\ref{3d vitro speed}) }\label{vitro-dR-dt}
   \end{minipage}\hfill
   \begin{minipage}{0.53\textwidth}
        \centering
        \includegraphics[width=1.0\linewidth]{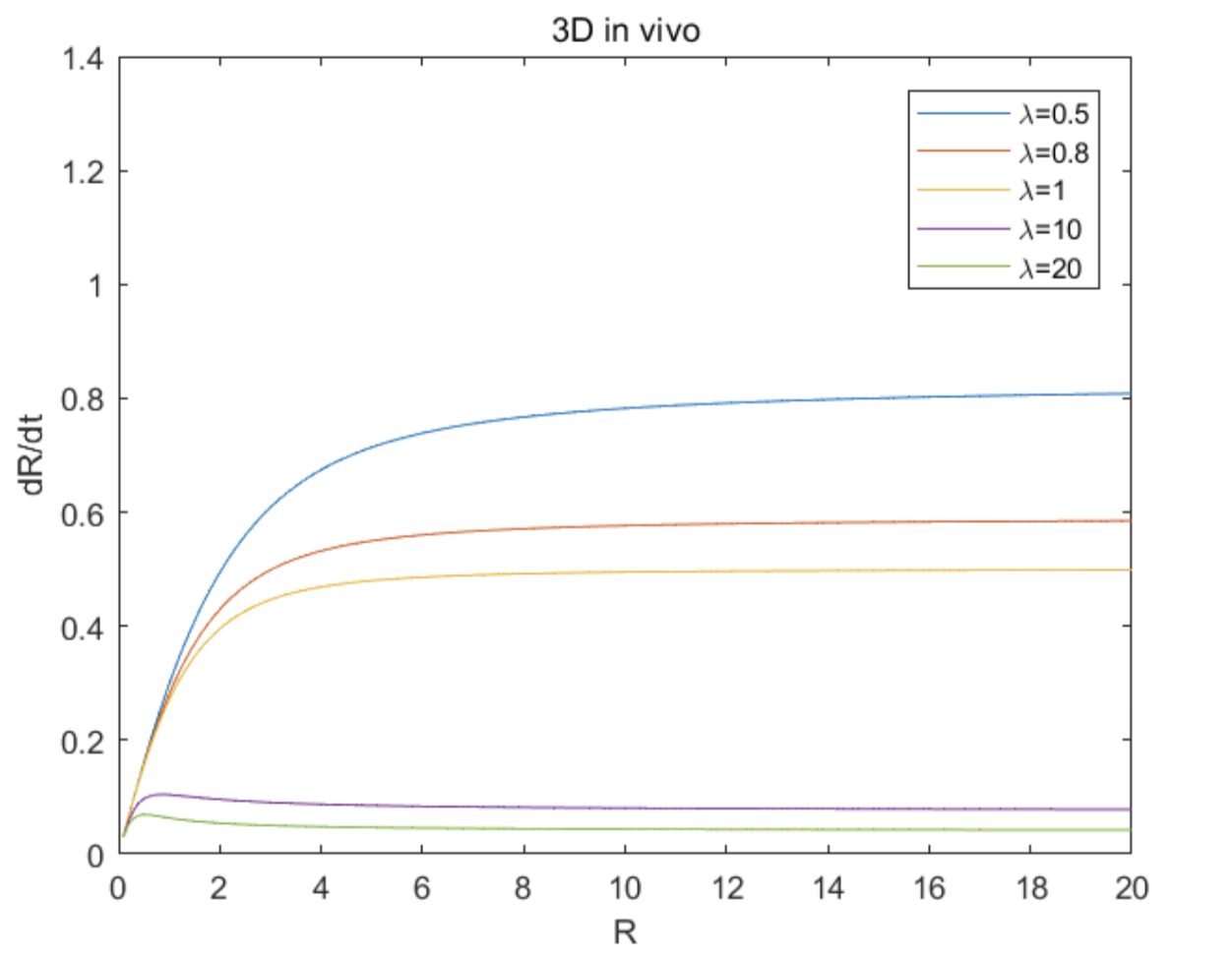}
        \caption{3D in vivo tumor evolution (\ref{3d speed in vivo})}\label{vivo-dR-dt}
   \end{minipage}
   \end{figure}

\begin{figure}[htp] 
        \centering
        \includegraphics[width=1.0\linewidth]{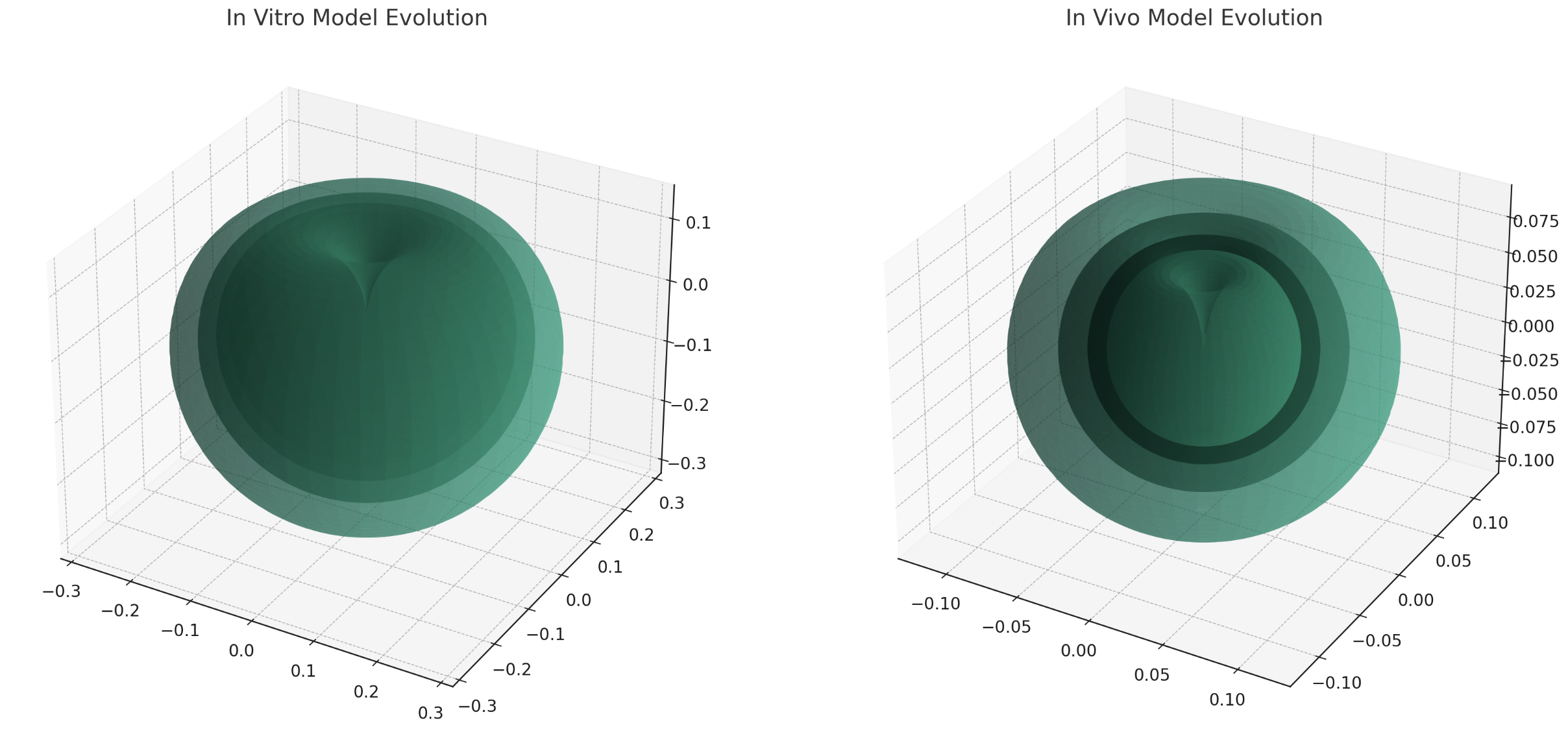}
        \caption{3D tumor evolution speed in vitro and in vivo (unperturbed) (\ref{3d vitro speed})(\ref{3d speed in vivo}) }\label{vitro-vivo-dR-dt}
   \end{figure}

We arrive at the Laplace's equation for $u=p+\frac{G_0}{\lambda}c$, done similarly to section 2, as shown in (\ref{useful combo}):
\begin{equation}
    \frac{d^2u}{dr^2}+\frac{2}{r}\frac{du}{dr}=0
\end{equation}
This is a Cauchy-Euler (CE) equation, and using the same approach as in (\ref{now 67 C-E}), we derive the general solution for pressure $p=B_1(R) - \frac{G_0}{\lambda}c$. Then we apply the boundary condition:
\begin{equation}
     \begin{aligned}
         p(r=R) = B_1(R) - \frac{G_0}{\lambda}c(R) = B_1(R) - \frac{G_0}{\lambda}c_B = 0 \\
         \Rightarrow B_1(R) = \frac{G_0 c_B }{\lambda}\frac{a_0(R)}{R}\sinh(\sqrt{\lambda} R)  =\frac{G_0 c_B}{\lambda}
     \end{aligned}
\end{equation}
Therefore, the solution for pressure $p$ for the in vivo model is 
\begin{equation} \label{p-3d-vivo}
    p(r,t) =  B_1(R) - \frac{G_0}{\lambda}c = \frac{G_0 c_B }{\lambda}\frac{a_0(R)}{R}(\sinh(\sqrt{\lambda} R)-\sinh(\sqrt{\lambda} r))
\end{equation}
where $a_0$ is shown in (\ref{c-3d-vivo}). \\ \\ 
\textbf{Step 3: Construct the boundary velocity} \\
Substitute the solution for $p$ we derived in (\ref{p-3d-vivo}) into 
\begin{equation}
    \frac{dR}{dt}=-\frac{\partial p}{\partial r} (R(t))
\end{equation}
We get the boundary velocity for the 3D spherically symmetric in vivo model is 
\begin{equation} \label{3d speed in vivo}
    \begin{aligned}
   \frac{dR}{dt}
   =\frac{G_0 c_B}{\lambda} \frac{R+1}{R^2}\frac{\cosh(\sqrt{\lambda}R) R \lambda - \sinh(\sqrt{\lambda}R)\sqrt{\lambda}}{\cosh(\sqrt{\lambda}R)\lambda+\sinh(\sqrt{\lambda}R)\sqrt{\lambda}}
    \end{aligned}
\end{equation}
Following equations (\ref{3d vitro speed}) and (\ref{3d speed in vivo}), we plot the comparison of the tumor evolution speed for in vitro and in vivo models as shown in Figures \ref{vitro-dR-dt} and \ref{vivo-dR-dt}, where we can observe that the tumor in the in vivo model is growing slower than the in vitro model. The 3D visualization of this is shown in Figure \ref{vitro-vivo-dR-dt}. Note that in Figure \ref{vitro-vivo-dR-dt}, the depressed part near $R=0$ may occur since the denominator has a higher order of $R$ in equations (\ref{3d vitro speed}) and (\ref{3d speed in vivo}). When the numerator increases at a slower rate than the denominator as $R$ gets smaller, this will cause the function to dip. When $R$ is zero or near zero, the denominator might approach zero faster than the numerator, leading to a dramatic increase which then quickly drops as $R$ increases and creates the depression part.
\subsection{\textbf{Stability of 3D spherically symmetric in the \textbf{in vitro} model}}
\leavevmode
For the spherically symmetric 3D case, we employ the spherical coordinates $(r,\theta, \phi)$, where we have 
\begin{equation}
    -\nabla^2 c(r,\theta, \phi, t) + \lambda c(r,\theta, \phi, t) = 0, \ \ \ \mbox{in}  \ \ \Tilde{D}
\end{equation}
and the boundary conditions are given as 
\begin{equation} \label{bcs in vitro c}
    \begin{aligned}
        \partial_rc_0(R(t),t)+c_1(R(t),t)=0 \\
    \end{aligned}
\end{equation}
\textbf{Step 1: Construct the concentration $c$.} \\
The spherical harmonic function  of degree $\ell$ and order $m$ is given by
\begin{equation} \label{Ylm}
    Y_{\ell}^{m}(\theta,\phi) = \sqrt{\frac{2\ell+1}{4 \pi}\frac{(\ell-m)!}{(\ell+m)!} }P_{\ell}^{m}(\cos\theta)e^{im\phi}
\end{equation}
where $\ell$ can take on any noninteger values, and for each value of $\ell$, $m$ can range from $-\ell$ to $\ell$.
$P_{\ell}^{m}(\cos\theta)$ is the associated Legendre polynomial. Note that $m=0$ mode corresponds to the unperturbed, spherically symmetric scenario. The stability analysis includes all modes, thereby generalizing to cases where $m=1,2,3,...$. This analysis consequently involves the associated Legendre polynomial $P_{\ell}^{m}(\cos\theta)$.
The 3D Laplacian in spherical coordinate $(r,\theta, \phi)$ is 
\begin{equation}
    \nabla^2 c = c_{rr} + \frac{2}{r}c_r + \frac{1}{r^2}\left[\frac{1}{\sin^2\theta}c_{\phi \phi} + \frac{1}{\sin \theta} (\sin\theta c_{\theta})_{\theta}   \right]
\end{equation}
After introducing a perturbation to the three-dimensional spherical tumor model, the series expansion for the concentration $c$, as established in \cite{feng2023tumor}, undergoes reduction. Due to the boundary conditions, all terms subsequent to the leading term in the series expansion for the concentration $c$ are identically zero, which considerably simplifies the expansion. Consequently, we are left with the following expression:
\begin{equation} \label{simplified exp}
    c_1(r,\theta,\phi,t)= c_1(r,t) Y_{\ell}^{m}(\theta,\phi)
\end{equation}
The 3D Laplacian is therefore
\begin{equation}
    \begin{aligned}    
        \nabla^2 c_1(r,\theta,\phi,t)=&(c_1)_{rr}Y_{\ell}^{m}(\theta,\phi)+\frac{2}{r}(c_1)_r Y_{\ell}^{m}(\theta,\phi)\\
        & + \frac{1}{r^2}\left[\frac{1}{\sin^2\theta}c_1 \frac{\partial^2 Y_{\ell}^{m}(\theta,\phi)}{\partial^2 \phi} + \frac{1}{\sin \theta} \partial_\theta\left(\sin \theta c_1 \frac{\partial Y_{\ell}^{m}(\theta,\phi)}{\partial \phi} \right)  \right]
    \end{aligned}
\end{equation}
Substituting equation (\ref{Ylm}) into the above 3D Laplacian in spherical coordinate, and since 
\begin{equation}
    \frac{d}{d\theta}\left(\sin\theta \frac{dc}{d\theta}\right) = -\left(\mu \sin\theta - \frac{\gamma}{\sin\theta} \right) c_1 Y_{\ell}^{m}(\theta,\phi)
\end{equation}
where $\mu=\ell(\ell+1)$ and $\gamma = m^2$.
We thus have
\begin{equation} \label{Laplace c perturbed}
    \begin{aligned}
        &\nabla^2 c_1(r,\theta,\phi) = (c_1)_{rr}Y_{\ell}^{m}(\theta,\phi)+\frac{2}{r}(c_1)_r Y_{\ell}^{m}(\theta,\phi) \\
        &+ \frac{1}{r^2 \sin^2 \theta} c_1 (-m^2) Y_{\ell}^{m}(\theta,\phi) + \frac{1}{r^2}\frac{1}{\sin\theta}\left(-\mu \sin\theta + \frac{\gamma}{\sin\theta}\right) c_1 Y_{\ell}^{m}(\theta,\phi)\\
        &=(r^2(c_1)_{rr}+2r(c_1)_r-\mu c_1)Y_{\ell}^m(\theta,\phi)
    \end{aligned}
\end{equation}
Therefore, plugging (\ref{Laplace c perturbed}) into (\ref{c real}) gives us 
\begin{equation}
    \begin{aligned}
     -\nabla^2 c+\lambda c = 0 \\
        \Rightarrow r^2(c_1)_{rr}+2r(c_1)_r-(\lambda r^2+\mu)c_1=0
    \end{aligned}
\end{equation}
where $\mu=\ell(\ell+1)$, and since the standard form for the modified spherical Bessel equation is
\begin{equation}
    x^2 \frac{d^2 y}{dx^2} + 2x\frac{dy}{dx} - (x^2+n(n+1))y=0
\end{equation}
We can convert
\begin{equation} \label{3d in vitro}
    r^2(c_1)_{rr}+2r(c_1)_r-(\lambda r^2+\mu)c_1=0
\end{equation}
to the standard form by change of variables, namely, $c(r)=y(z)/\sqrt{z}$ with $z=\sqrt{\lambda}r$. 
Therefore, the solution to Equation \eqref{3d in vitro} is
\begin{equation}
    \begin{aligned}
        y(z) = C_1 I_{w}(z)+C_2 K_{w}(z) \\
    \end{aligned}
\end{equation}
where $w=\ell+\frac{1}{2}$. Since the solution is bounded at $z=0$, which sets $c_2=0$, then 
\begin{equation}
    c(r)= \frac{C_1 I_{\ell+\frac{1}{2}}(z)}{\sqrt{\sqrt{\lambda}r}}
\end{equation}
Since we have the boundary conditions:
\begin{equation} \label{3d vitro concentration}
    \begin{aligned}
        \partial_r c_0(R(t),t)+c_1(R(t),t)=0 \\
        c_1(R(t),t)=0
    \end{aligned}
\end{equation}
Thus from (\ref{3d vitro concentration}):
\begin{eqnarray} \label{c11 vitro}
     c_1(R(t),t)
     &=& -\frac{c_B \sqrt{\lambda} \cosh(\sqrt{\lambda}r)}{\sinh(\sqrt{\lambda}R)} \nonumber
\end{eqnarray}
By using the boundary conditions in (\ref{bcs in vitro c}), we have
\begin{equation}
     -\frac{c_B \sqrt{\lambda} \cosh(\sqrt{\lambda}r)}{\sinh(\sqrt{\lambda}R)} = C_1 \frac{I_{\ell+\frac{1}{2}}(\sqrt{\lambda}r)}{\sqrt{\sqrt{\lambda}r}} 
\end{equation}
Therefore, for the in vitro model, we have the concentration inside the tumor: 
\begin{equation}
    c_1(r,t) = c_B a_1(t)\frac{I_{\ell+\frac{1}{2}}(\sqrt{\lambda}r)}{\sqrt{\sqrt{\lambda}r}}
\end{equation}
and the leading-order term $c_0(r,t)$ for this case is given by (\ref{3d vitro concentration}). We have 
\begin{equation}
    a_1(R(t))=-\frac{\sqrt{\lambda}\cosh(\sqrt{\lambda }R) \sqrt{\sqrt{\lambda}R}}{\sinh(\sqrt{\lambda R})I_{\ell+\frac{1}{2}}(\sqrt{\lambda}R)}
\end{equation}
The concentration outside the tumor is the constant $c_B$. \\ \\
\textbf{Step 2: Construct pressure $p_1$.}\\
The PDE for the pressure $p$ is 
\begin{equation}
    -\nabla^2 p_1(r,\theta,\phi,t)=G_0 c_1(r,\theta,\phi,t), \ \ \mbox{in} \ \ \Tilde{D}(t) 
\end{equation}
with boundary condition
\begin{equation} \label{3d p perturb bcs}
    \begin{aligned}
        \partial_r p_0(R(t),t)+p_1(R(t),t)=0, 
    \end{aligned}
\end{equation}
Since $-\nabla^2 c_1(r, \theta, \phi, t)+\lambda c_1(r, \theta, \phi, t)=0$, meaning that $c_1=\frac{1}{\lambda}\nabla^2 c_1$, therefore
\begin{equation}
    \begin{aligned}
        -\nabla^2 p_1(r, \theta, \phi, t)=Gc_1(r, \theta, \phi, t) \\
        -\nabla^2 p_1 - \frac{G}{\lambda}\nabla^2 c_1 = 0 
    \end{aligned}
\end{equation}
Then we substitute the Laplacian in 3D spherical coordinate for pressure $p$, combining the pressure and concentration equations as in (\ref{useful combo}), which gives us the Laplace's equation for $u_1=p_1+\frac{G}{\lambda}c_1$. Since we have the Laplace operator after perturbation as in (\ref{Laplace c perturbed}), we have:
\begin{equation}
    r^2 \frac{d^2(u_1)}{dr^2} + 2r \frac{du_1}{dr}-\mu u_1=0
\end{equation}
This is in Cauchy-Euler form, and has the general solution 
\begin{equation}
    u=H_1 r^{\ell}+H_2 r^{-\ell}
\end{equation}
Since the solution is bounded, hence we get
\begin{eqnarray}
    u=H_1(R)r^{\ell}
\end{eqnarray}
By applying the boundary conditions in (\ref{3d p perturb bcs}), and substitute $$p_0=-\frac{G_0}{\lambda}\frac{1}{r}\frac{c_B R \sinh(\sqrt{\lambda}r)}{\sinh(\sqrt{\lambda}R}+\frac{G_0}{\lambda}c_B$$ as in (\ref{p-3d-vitro}). Therefore,
\begin{equation}
    p_1(R(t),t)=-\partial_r p_0(R(t),t)=-\frac{R G_0 c_B \sinh(\sqrt{\lambda}R)}{R^2 \lambda \sinh(\sqrt{\lambda}R)}+\frac{R G_0 c_B \cosh(\sqrt{\lambda}R)}{R \sqrt{\lambda}\sinh(\sqrt{\lambda}R)}
\end{equation}
Since $p_1=u_1-\frac{G_0 c_B}{\lambda}$, we have
\begin{equation} \label{p11=u11-}
    p_1(R,t) = \frac{G_0 c_B}{\lambda} \left(-\frac{1}{R} + \frac{\sqrt{\lambda} \cosh(\sqrt{\lambda}R)}{\sinh(\sqrt{\lambda}R)} \right)
\end{equation}
Substitute (\ref{c11 vitro}) for $c_1(R(t),t)$ in the above equation and solve for $H_1(R)$, we get,
\begin{equation} \label{B_1}
    \begin{aligned}
        H_1(R) &= R^{-\ell} \left[ \left( \frac{R G_0 c_B \sinh(\sqrt{\lambda}r)}{r^2 \lambda \sinh(\sqrt{\lambda}R)}\!-\!\frac{R G_0 c_B \cosh(\sqrt{\lambda}r)}{r \sqrt{\lambda}\sinh(\sqrt{\lambda}R)}\right) \!-\! \frac{G_0}{\lambda}\left( \frac{c_B \sqrt{\lambda}\cosh(\sqrt{\lambda}R)}{\sinh(\sqrt{\lambda}R)}\right)  \right] \\
        &\quad \Rightarrow H_1(R) = -R^{-\ell-1} \left( \frac{G_0 c_B \sinh(\sqrt{\lambda}R)}{\lambda \sinh(\sqrt{\lambda}R)} \right)
        =-R^{-\ell-1}\left( \frac{G_0 c_B}{\lambda} \right)
    \end{aligned}
\end{equation}
We then substitute $(\ref{B_1})$ into (\ref{p11=u11-}), and we get the solution for $p_1$,
\begin{equation} \label{p11 final solution}
 p_1(r,t)=-\frac{G_0 c_B}{\lambda} \left(     -\frac{r^{\ell}}{R^{\ell+1} }  +  \frac{\sqrt{\lambda}\cosh(\sqrt{\lambda}R) \sqrt{R}}{\sinh(\sqrt{\lambda}R) I_{\ell+\frac{1}{2}}(\sqrt{\lambda}R)} \frac{I_{\ell+\frac{1}{2}}(\sqrt{\lambda}r)}{\sqrt{r}} \right)
\end{equation}
\textbf{Step 3: Construct $\delta^{-1} \frac{d \delta}{dt}$.} \\
From the previous perturbation analysis we have the evolution function
\begin{equation}
    \delta^{-1}\frac{d \delta}{dt} = -(\partial^2_r p_0(R,t)+\partial_r p_1(R,t) + O(\delta))
\end{equation}
We substitute $p_0(R)=-\frac{G_0}{\lambda}\frac{1}{r}\frac{c_B R \sinh(\sqrt{\lambda}R)}{\sinh(\sqrt{\lambda}R)}+\frac{G_0 c_B}{\lambda}$
and $p_1$ as in (\ref{p11 final solution}), so we get 
\begin{equation} \label{3d vitro evol}
    \begin{aligned}      
        \delta^{-1}\frac{d \delta}{dt} &=  \frac{ R c_B G_0  I_{\ell+\frac{1}{2}}(\sqrt{\lambda}R) R^{-1-\ell}\sinh(\sqrt{\lambda}R) R^{\ell+\frac{5}{2}}\ell \sqrt{\lambda}}{\lambda^{3/2}R^{7/2}\sinh(\sqrt{\lambda}R) I_{\ell+\frac{1}{2}}(\sqrt{\lambda}R) }\\
        & +\frac{R c_B G_0  I_{\ell+\frac{1}{2}}(\sqrt{\lambda}R) (2\sqrt{\lambda}\sqrt{R}+\lambda^{3/2} R^{5/2})\sinh(\sqrt{\lambda}R)
        }{\lambda^{3/2}R^{7/2}\sinh(\sqrt{\lambda}R) I_{\ell+\frac{1}{2}}(\sqrt{\lambda}R)}  \\
        & -\frac{2R c_B G_0  I_{\ell+\frac{1}{2}}(\sqrt{\lambda}R) \cosh(\sqrt{\lambda}R)\lambda R^{3/2}}{\lambda^{3/2}R^{7/2}\sinh(\sqrt{\lambda}R) I_{\ell+\frac{1}{2}}(\sqrt{\lambda}R)}  \\
        &  -\frac{ (-\lambda(\ell+1)I_{\ell+\frac{1}{2}}(\sqrt{\lambda}R)+I_{\ell-\frac{1}{2}}(\sqrt{\lambda}R) \lambda^{3/2}R) R^{5/2} \cosh(\sqrt{\lambda}R) 
        } {2\lambda^{3/2}R^{7/2}\sinh(\sqrt{\lambda}R) I_{\ell+\frac{1}{2}}(\sqrt{\lambda}R) } +O(\delta) \\
        & \defeq M_1(\lambda,\ell,R)+O(\delta)
    \end{aligned}
\end{equation}
\textbf{Step 4: Boundary stability analysis.} \\
The goal is to determine the sign of the evolution function $\delta^{-1} \frac{d \delta}{dt}$, and then determine the stability of the tumor boundary. By substituting the asymptotic expansions for the modified Bessel functions and the hyperbolic functions as in Appendix A, we arrive at (\ref{M1}) in section 5, which are summarized as the Corollary 5.5 and indicate the negative sign of the evolution function in 3D in vitro scenario. From Figures \ref{3d vitro-1},\ref{3d vitro-2} below, we can also see that the evolution function $\delta^{-1} \frac{d \delta}{dt} <0$ is always true, meaning that in the 3D in vitro model, the perturbation will degenerate and the tumor boundary will always become stable.
\begin{figure}[htp] 
    \begin{minipage}{0.53\textwidth}
        \centering
        \includegraphics[width=1.0\linewidth]{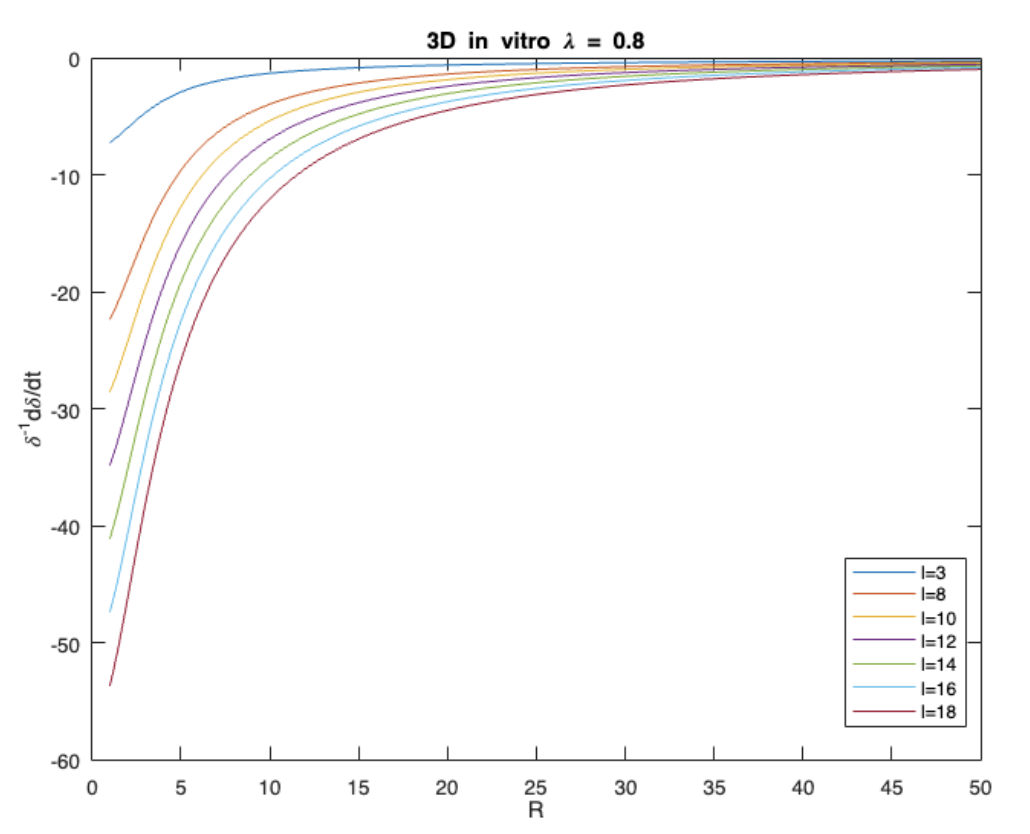}
        \caption{3D in vitro $\lambda=0.8,\ell=3,8,12,16,18$} \label{3d vitro-1}
   \end{minipage}\hfill
   \begin{minipage}{0.53\textwidth}
        \centering
        \includegraphics[width=1.0\linewidth]{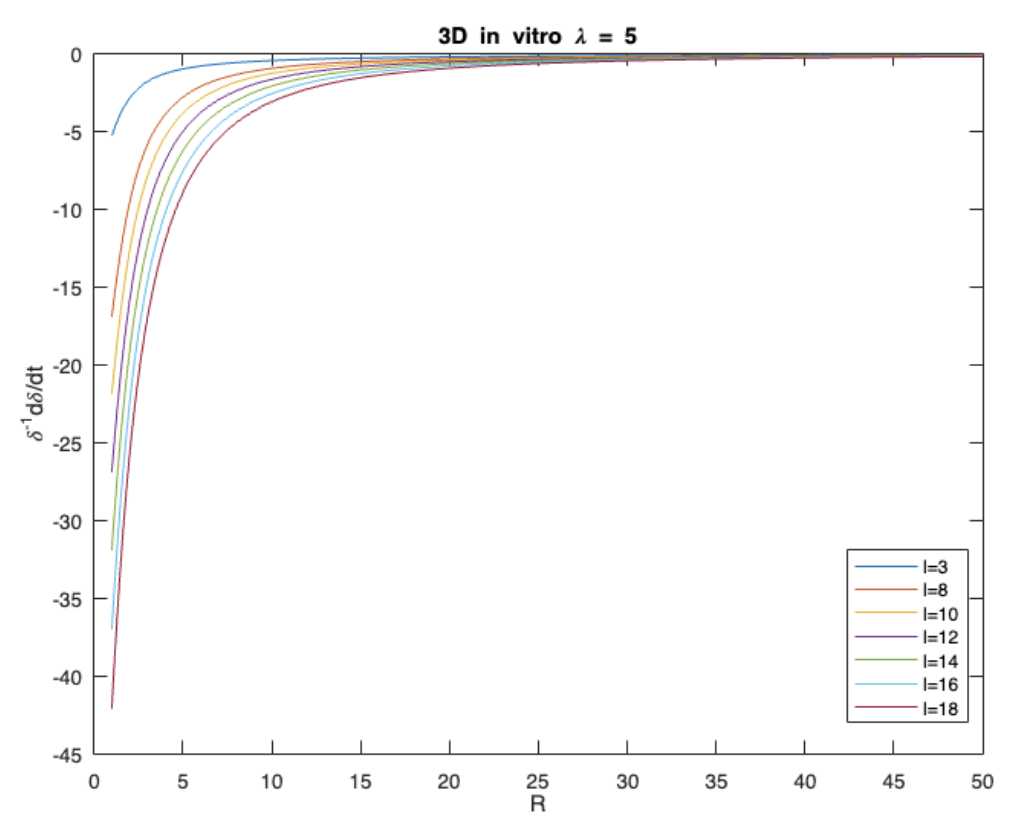}
        \caption{3D in vitro $\lambda=5,\ell=3,8,12,16,18$}\label{3d vitro-2}
   \end{minipage}
   \end{figure}

\subsection{\textbf{Stability of 3D spherically symmetric in the \textbf{in vivo} model}}
The PDEs for the 3D in vivo problem is 
\begin{equation} \label{govern 3d in vivo}
    \begin{aligned}
        -\nabla^2 c_1^{(i)}(r, \theta, t)+\lambda c_1^{(i)}(r, \theta, t) = 0, \\
        -\nabla^2 c_1^{(o)}(r, \theta, t)+ c_1^{(o)}(r, \theta, t) = 0 \\
    \end{aligned}
\end{equation}
with the solution of concentration ($c_1^{(i)}(0,t)$) bounded everywhere, the boundary conditions (BCs) and initial conditions (ICs):
\begin{equation} \label{3d vivo bcs perturb}
    \begin{aligned}
        c_1^{(o)}(\infty,t)=0,  \\
        c_1^{(i)}(R(t),t)=c_1^{(o)}(R(t),t), \\
        \partial_rc_1^{(i)}(R(t),t)=\partial_r c_1^{(o)}(R(t),t), \\
        \partial_r^2 c_0^{(i)}(R(t),t)+\partial_rc_1^{(i)}(R(t),t)=\partial_r^2 c_0^{(o)}(R(t),t)+\partial_rc_1^{(o)}(R(t),t) 
    \end{aligned}
\end{equation}
\textbf{Step 1: Construct the concentration $c$.} \\
For the in vivo model, the idea is to plug the simplified expansion as in (\ref{simplified exp}): 
$$c_1(r,\theta,t)=   c_1(r,t) Y_{\ell}^{m}(\theta,\phi) $$ 
into (\ref{govern 3d in vivo}), where $\mu=\ell(\ell+1), \gamma=m^2$. \\

The construction for the concentration inside the tumor $c^{(i)}(r,t)$ is exactly the same as that in the in vitro model when solving $c(r,t)$, and thus we skip the construction of $c^{(i)}(r,t)$ and arrive at the general solution 
\begin{equation}  \label{3d vivo inside c}
    c^{(i)}(r,t)=C_1 \frac{I_{\ell+\frac{1}{2}}(\sqrt{\lambda}r)}{\sqrt{\sqrt{\lambda}r}} 
\end{equation}
Then we consider the concentration outside the tumor $c^{(o)}(r,t)$ governed by the second equation in (\ref{govern 3d in vivo}) and since we have $\nabla^2 c_1(r,\theta,\phi)=r^2(c1)_{rr}+2r(c_2)_r- \mu c_1$ as derived in (\ref{Laplace c perturbed}), we get:
\begin{eqnarray}
    \Rightarrow \left[r^2(c_1)_{rr}+2r(c_1)_r-\left(\frac{{m}^2}{r^2\sin^2\theta} + \frac{\mu}{r^2} - \frac{\gamma}{r^2 \sin^2\theta}+1\right)c_1\right]Y_{\ell}^m(\theta,\phi)&=&0 \\ \nonumber
    \Rightarrow r^2(c_1)_{rr}+2r(c_1)_r-(r^2+\mu)c_1&=&0
\end{eqnarray}
where $\mu=\ell(\ell+1)$, and since the standard form for the modified spherical Bessel function is
\begin{equation}
    x^2 \frac{d^2 y}{dx^2} + 2x\frac{dy}{dx} - (x^2+n(n+1))y=0
\end{equation}

We can convert
\begin{equation} 
    r^2(c_1)_{rr}+2r(c_1)_r-(r^2+\mu)c_1=0
\end{equation}
to the standard form by change of variables, namely, $c(r)=y(r)/\sqrt{r}$. Therefore, the solution to the governing equation for concentration outside the tumor is
\begin{equation}
    c^{(o)}(r,t)=C_3 \frac{I_{m}(r)}{\sqrt{r}}+C_4 \frac{K_{m}(r)}{\sqrt{r}}
\end{equation}
where $m=\ell+\frac{1}{2}$ as before. Since $I_m$ blows up when $r \rightarrow \infty$, $C_3=0$, and the general solution becomes 
\begin{equation} \label{3d vivo outside c}
    c^{(o)}(r,t)=C_4 \frac{K_{m}(r)}{\sqrt{r}}
\end{equation}

Then we apply the boundary conditions as in (\ref{3d vivo bcs perturb}) to get the constants $C_1, C_4$ and write in a consistent form. Hence, for any $k \in \mathbb{N^+}$ we have 
\begin{eqnarray} \label{ci,co 3d vivo-1}
    c_1^{(i)}(r,t)&=&c_B a_1(t) \frac{I_{\ell+\frac{1}{2}}(\sqrt{\lambda}r)}{\sqrt{\sqrt{\lambda}r}}, \ \ \mbox{for} \ \ r \leq R(t), \\
    c_1^{(o)}(r,t)&=&c_B b_1(t) \frac{K_{\ell+\frac{1}{2}}(r)}{\sqrt{r}}, \ \ \mbox{for} \ \ r \geq R(t).
\end{eqnarray}
Recall that the leading-order terms $c_0^{(i)}(\xi)$ and $c_0^{(o)}(\xi)$ are given by (\ref{also 3d c}) and (\ref{c-3d-vivo}). Then applying the boundary conditions (\ref{3d vivo bcs perturb}) yields:
\begin{equation}\label{a11fv}
    \begin{aligned}
        a_1(t) &=\frac{(((R+1)\lambda^{3/2}-\sqrt{\lambda})\sinh{(\sqrt{\lambda}R)}+\cosh{(\sqrt{\lambda}R)}R\lambda) (\sqrt{\lambda}R)^{3/2}K_{\ell+\frac{1}{2}}(R)}{\lambda R C(R)} \\
    \end{aligned}    
\end{equation}
\begin{equation} \label{b11fv}
    \begin{aligned}
        b_1(t) &=\frac{2(((R+1)\lambda^{3/2}-\sqrt{\lambda})\sinh{(\sqrt{\lambda}R)}+\cosh{(\sqrt{\lambda}R)}R\lambda)I_{\ell+\frac{1}{2}}(\sqrt{\lambda}R)}{2C(R)}\\
    \end{aligned}
\end{equation}
where,
\begin{equation} 
    \begin{aligned}
        C(R)&=(\sqrt{\lambda}\cosh{(\sqrt{\lambda}R)} + \sinh{(\sqrt{\lambda}R)})\left( K_{\ell-\frac{1}{2}}(R)R-\frac{ K_{\ell+\frac{1}{2}}(R)}{2}   \right)I_{\ell+\frac{1}{2}}(\sqrt{\lambda}R)\\
        &+RI_{\ell-\frac{1}{2}}(\sqrt{\lambda}R) K_{\ell+\frac{1}{2}}(R) (\cosh{(\sqrt{\lambda}R)}\lambda + \sinh{(\sqrt{\lambda}R)}\sqrt{\lambda})\\
    \end{aligned}
\end{equation}\\
\textbf{Step 2: Construct pressure $p_1$.}\\
By now, $c_1(r,t)$ is determined, therefore, $c_1(r,\theta,\phi)$ is also determined. Then, by solving the equation for pressure $p$:
\begin{equation}
    -\nabla^2 p_1(r,\theta,\phi,t)=G_0 c_1(r,\theta,\phi,t), \ \ \mbox{in} \ \ \Tilde{D}(t) 
\end{equation}
with boundary condition
\begin{equation} \label{3d p perturb bcs_2}
    \begin{aligned}
         \partial_r p_0(R(t),t)+p_1(R(t),t)=0, 
    \end{aligned}
\end{equation}
and expansion $p_1(r,\theta,t)=   p_1(r,t) Y_{\ell}^m(\theta,\phi) $ and $p_0$ given by (\ref{p-3d-vivo}).
we can solve for $p_1$ . 

Specifically, since $-\nabla^2 c_1(r, \theta, \phi, t)+\lambda c_1(r, \theta, \phi, t)=0$, meaning that $c_1=\frac{1}{\lambda}\nabla^2 c_1$, therefore
\begin{equation}
    \begin{aligned}
        -\nabla^2 p_1(r, \theta, \phi, t)=G_0c_1(r, \theta, \phi, t) \\
        -\nabla^2 p_1 - \frac{G}{\lambda}\nabla^2 c_1 = 0 
    \end{aligned}
\end{equation}
Then we substitute the Laplacian in 3D spherical coordinate for pressure $p$, combining the pressure and concentration equations as in (\ref{useful combo}), which gives us the Laplace's equation for $u_1=p_1+\frac{G}{\lambda}c_1$. Since we have the Laplace operator after perturbation as in (\ref{Laplace c perturbed}), we have:
\begin{equation}
    r^2 \frac{d^2 u_1}{dr^2} + 2r \frac{du_1}{dr}-\mu u_1=0
\end{equation}
This is in Cauchy-Euler form, and has the general solution $u_1=D_1(R) r^{\ell}+D_2(R) r^{\ell}$, since the solution is bounded at the origin, $D_2(R)=0$, and we have $u_1=D_1(R) r^{\ell}$.\\
By applying the boundary conditions and substitute $$p_0=  \frac{G_0 c_B }{\lambda}\frac{a_0(R)}{R}(\sinh(\sqrt{\lambda} R)-\sinh(\sqrt{\lambda} r))$$ as in (\ref{p-3d-vivo}) and $a_0$ given in (\ref{c-3d-vivo}). Therefore,
\begin{equation}
    \begin{aligned} \label{p11new}
         p_1(R(t),t)&=-\partial_r p_0(R(t),t) \\
&=\frac{G_0 c_B}{\lambda} \frac{R+1}{R^2}\frac{\cosh(\sqrt{\lambda}R) R \lambda - \sinh(\sqrt{\lambda}R)\sqrt{\lambda}}{\cosh(\sqrt{\lambda}R)\lambda+\sinh(\sqrt{\lambda}R)\sqrt{\lambda}}
    \end{aligned}
\end{equation}
and since $p_1(R,t) = u_1(R)-G_0 c_1/\lambda$, we have
\begin{equation} 
    \begin{aligned} \label{p11-vivo}
        p_1(R,t)=D_1(R) R^{\ell}-\frac{G_0}{\lambda}c_1(R) \\
    \end{aligned} 
\end{equation}
We combine (\ref{p11new}) and (\ref{p11-vivo}), and substitute (\ref{ci,co 3d vivo-1}) for $c_1^{(i),1}(R(t),t)$ in the above equation to solve for constant $D_1(R)$, we get 
\begin{equation}\label{G_1-vivo}
\begin{split}
     D_1(R) = \frac{G_0 c_B}{\lambda} \left(   \frac{R+1}{R^{2+\ell}} \frac{\cosh(\sqrt{\lambda}R)R\lambda-\sinh(\sqrt{\lambda}R)\sqrt{\lambda}}{\cosh(\sqrt{\lambda}R)\lambda+\sinh(\sqrt{\lambda}R)\sqrt{\lambda}}+ \frac{a_1}{R^{\ell}}\frac{I_{\ell+\frac{1}{2}}(\sqrt{\lambda}R)}{ \sqrt{\sqrt{\lambda}R}}         \right)
\end{split}
\end{equation}
where $b_1$ is shown in (\ref{b11fv}).\\ \\
We then substitute $(\ref{G_1-vivo})$ and $(\ref{b11fv})$ into (\ref{p11-vivo}), and we get the solution for $p_1$ for the in vivo model
\begin{equation} \label{p11 final sol}
    \begin{aligned}
        p_{1}(t)= \frac{G_0 c_B}{\lambda} \left[ \left(   \frac{R+1}{R^{2+\ell}} \frac{\cosh(\sqrt{\lambda}R)R\lambda-\sinh(\sqrt{\lambda}R)\sqrt{\lambda}}{\cosh(\sqrt{\lambda}R)\lambda+\sinh(\sqrt{\lambda}R)\sqrt{\lambda}}+ \frac{a_1}{R^{\ell}}\frac{I_{\ell+\frac{1}{2}}(\sqrt{\lambda}R)}{ \sqrt{\sqrt{\lambda}R}}         \right) r^{\ell} \right. \\
        \left. - a_1 \frac{I_{\ell+\frac{1}{2}}(\sqrt{\lambda}R) }{\sqrt{\sqrt{\lambda}R}} \right]
    \end{aligned}
\end{equation}
where $a_1$ is shown in (\ref{a11fv}).\\\\
\textbf{Step 3: Construct $\delta^{-1} \frac{d \delta}{dt}$.} \\
Since we have the expression for the evolution function $\delta^{-1}\frac{d \delta}{dt}$:
\begin{equation} \label{del}
    \begin{aligned}
        \delta^{-1}\frac{d \delta}{dt} = -(\partial_r^2p_0(R,t)+\partial_r p_1(R,t))
    \end{aligned}
\end{equation}
with
\begin{equation}
    \begin{aligned}
        p_0(r,t) &= \frac{-G_0 c_B}{\lambda} \left( \frac{a_0(R)\sinh(\sqrt{\lambda}r)}{r}+\frac{a_0(R)\sinh(\sqrt{\lambda}R) }{R} \right)
    \end{aligned}
\end{equation}
where $a_0$ is constructed in (\ref{c-3d-vivo}). $p_1$ is constructed as in (\ref{p11 final sol}). \\
We substitute them into (\ref{del}), which gives us
\begin{equation} \label{3d-vivo-delta}
    \delta^{-1}\frac{d\delta}{dt} \defeq \frac{T_1}{T_2}+O(\delta) \defeq M_2(\lambda,\ell,R)+O(\delta)
\end{equation}
where,
\begin{equation} \label{T1}
    \begin{aligned}
        T_1& = - G_0 c_B( ((-2R( (3R+2)\lambda^{3/2}+\sqrt{\lambda}(R-2) )(\ell+1/2) K_{\ell+\frac{1}{2}}(R) \\
        &-2(R(R-\ell/2-1) \lambda^{3/2}-\sqrt{\lambda}(\ell+2)(R-2)/2)(R+1)K_{\ell-\frac{1}{2}}(R))\cdot  \\
        &\cosh{(\sqrt{\lambda}R)}^3-2\sinh{(\sqrt{\lambda}R)} ( R(\ell+1/2)( (\lambda^2+3\lambda) R+\lambda^2-1 ) K_{\ell+\frac{1}{2}}(R) \\
        &+((\lambda^2+\lambda)R^2-2\lambda(\ell+2)R+(\lambda+1)(\ell+2))(R+1)K_{\ell-\frac{1}{2}}(R)/2)\cdot\\ 
        & \cosh{(\sqrt{\lambda}R)}^2+(2R((2R+2)\lambda^{3/2}+\sqrt{\lambda}(R-2))(\ell+1/2)K_{\ell+\frac{1}{2}}(R) \\
        &+2(R^2\lambda^{3/2}-\sqrt{\lambda}(\ell+2)(R-2)/2)(R+1)K_{\ell-\frac{1}{2}}(R))\cosh{(\sqrt{\lambda}R)} \\
        & +2\sinh{(\sqrt{\lambda}R)}(R(\ell+1/2)(\lambda R+\lambda-1)K_{\ell+\frac{1}{2}}(R)\\
        &+K_{\ell-\frac{1}{2}}(R)(R+1)(R^2\lambda+\ell+2)/2))I_{\ell+\frac{1}{2}}(\sqrt{\lambda}R)\\
        &+I_{\ell-\frac{1}{2}}(\sqrt{\lambda}R)(R^2+(\ell+2)R+\ell+2)(\lambda(-2+(\lambda+1)R)\cosh{(\sqrt{\lambda}R)}^3\\
        &+2\sinh{(\sqrt{\lambda}R)}((R-1/2)\lambda^{3/2}-\sqrt{\lambda}/2)\cosh{(\sqrt{\lambda}R)}^2\\
        &-\lambda(R-2)\cosh{(\sqrt{\lambda}R)}+\sinh{(\sqrt{\lambda}R)}\sqrt{\lambda})K_{\ell+\frac{1}{2}}(R) \\
    \end{aligned}
\end{equation}
and,
\begin{equation} \label{T2}
    \begin{aligned}  
        T_2&= \sqrt{\lambda}R^3 K_{\ell+\frac{1}{2}}(R) (\cosh{(\sqrt{\lambda}R)}\lambda+ \sinh{(\sqrt{\lambda}R)}\sqrt{\lambda} )I_{\ell-\frac{1}{2}}(\sqrt{\lambda}R)\\
        &+I_{\ell+\frac{1}{2}}(\sqrt{\lambda}R) K_{\ell-\frac{1}{2}}(R)(\sqrt{\lambda} \cosh{(\sqrt{\lambda}R)} + \sinh{(\sqrt{\lambda}R)}))\cdot\\
        & (\cosh{(\sqrt{\lambda}R)}\lambda+ \sinh{(\sqrt{\lambda}R)}\sqrt{\lambda} )   (\sqrt{\lambda} \cosh{(\sqrt{\lambda}R)} + \sinh{(\sqrt{\lambda}R)})\\
    \end{aligned}
\end{equation}
\textbf{Step 4: Boundary stability analysis.} \\
We have evolution function as in (\ref{3d-vivo-delta})-(\ref{T2}). To determine the sign of the evolution function, we can substitute the asymptotic expansions of the modified Bessel functions and the hyperbolic trigonometric functions when $r\to \infty$ and $r \to 0$ as shown in Appendix.\ref{secAAA}. Specifically,
in our analysis, we first substitute the given asymptotic expansions into the evolution function. This step is followed by simplifying the expression and arranging the terms in descending order of powers of \( R \). For the case when \( r \to \infty \), we focus on the highest power term of \( R \) due to its dominance in the evolution function for large \( R \) values. We extract its coefficient as the key term. Conversely, when \( r \to 0 \), the lowest power term of \( R \) becomes significant. Here, we isolate this term from the evolution function after the substitution of the relevant asymptotic expansions and retrieve its coefficient. This analysis clarifies the evolution function's behavior near zero and infinity, pivotal for determining tumor boundary stability; a negative evolution function denotes stability, whereas a positive evolution function indicates instability. Results of this will be shown in section 5.1.
\begin{figure}[htp] 
   \begin{minipage}{0.52\textwidth}
        \centering
        \includegraphics[width=1.0\linewidth]{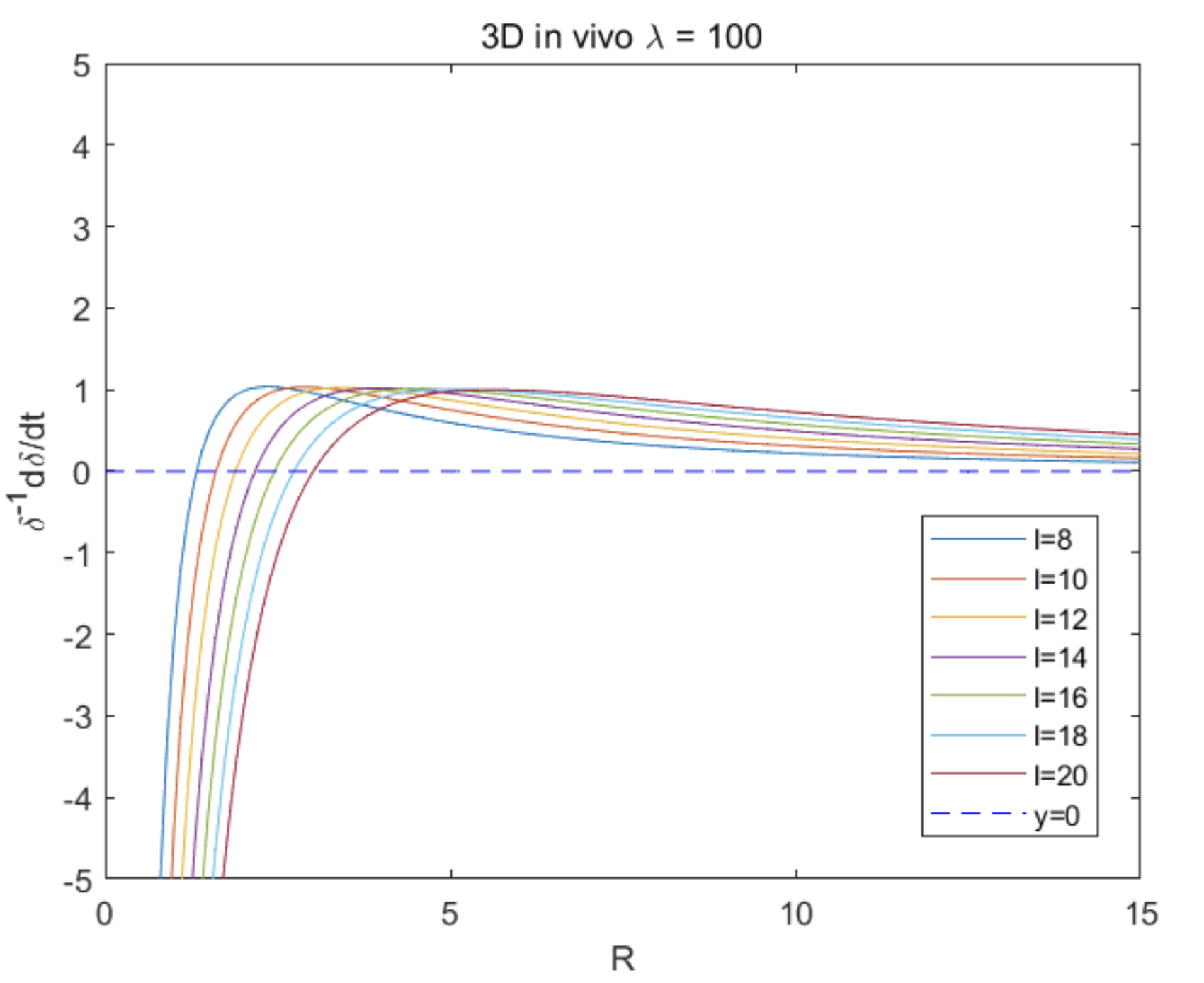}
        \captionsetup{justification=centering}
        \caption{3D in vivo, $G_0=1, c_B=100, \lambda=100, \ell=8,10,12,14,16,18,20, R \in [0,15] $}\label{12}
   \end{minipage}\hfill
   \begin{minipage}{0.52\textwidth}
        \centering
        \includegraphics[width=1.0\linewidth]{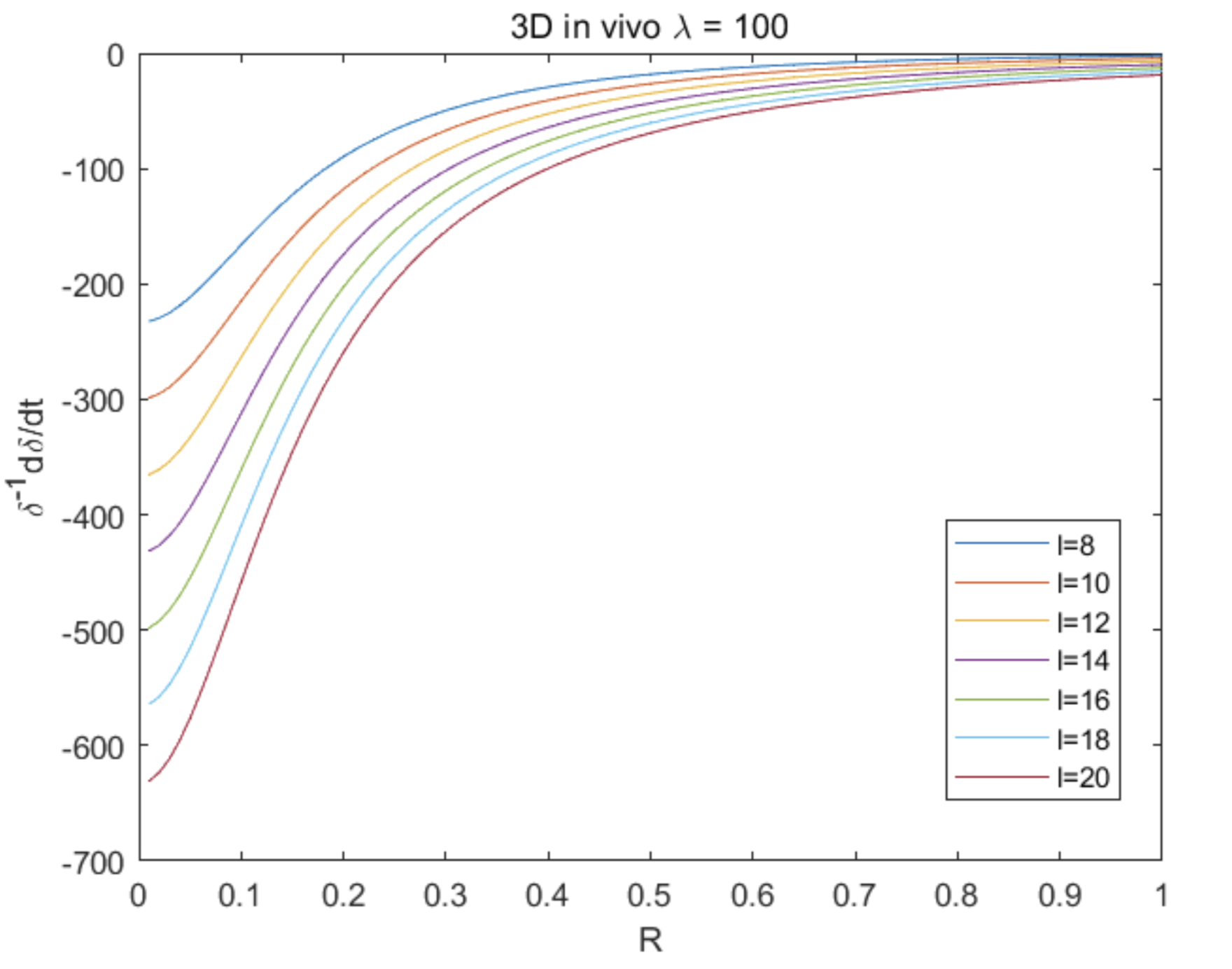}
        \caption{3D in vivo, $G_0=1,c_B=100, \lambda=100, \ell=8,10,12,14,16,18,20, R \in [0,1]$}\label{13}
   \end{minipage}
   
   \begin{minipage}{0.52\textwidth}
        \centering
        \includegraphics[width=1.0\linewidth]{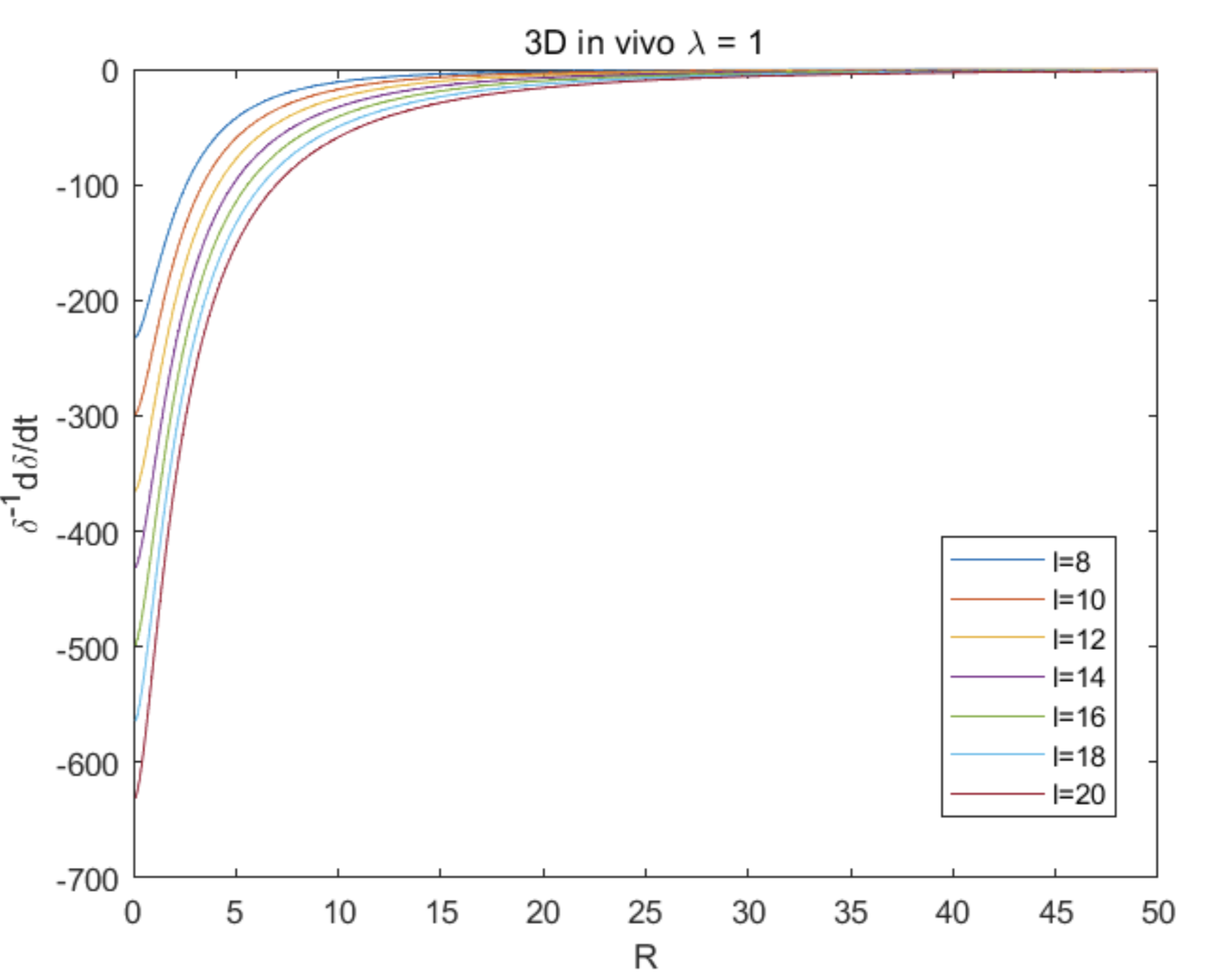}
        \captionsetup{justification=centering}
        \caption{3D in vivo, $G_0=1, c_B=100, \lambda=1.0, \ell=8,10,12,14,16,18,20, R \in [0,50]$} \label{14}
   \end{minipage}\hfill
   \begin{minipage}{0.52\textwidth}
        \centering
        \includegraphics[width=1.0\linewidth]{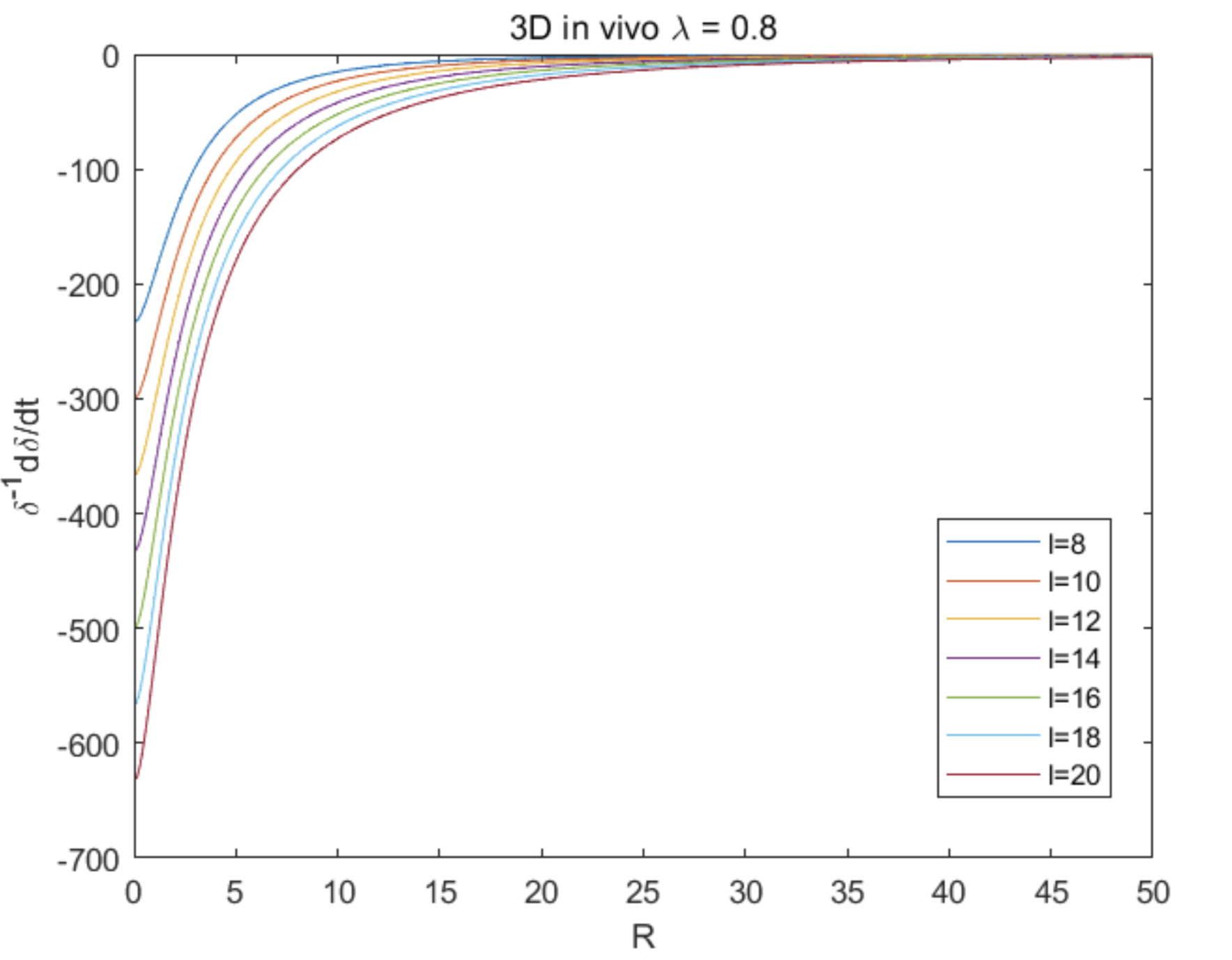}
        \caption{3D in vivo, $G_0=1, c_B=100, \lambda=0.8, \ell=8,10,12,14,16,18,20, R\in[0,50]$} \label{15}
   \end{minipage}
   \end{figure}
  
\section{Results and Discussion}\label{sec12}
\subsection{Results for 3D tumor model}
\begin{figure}[htp] 
    \begin{minipage}{0.53\textwidth}        \centering
        \includegraphics[width=1.0\linewidth]{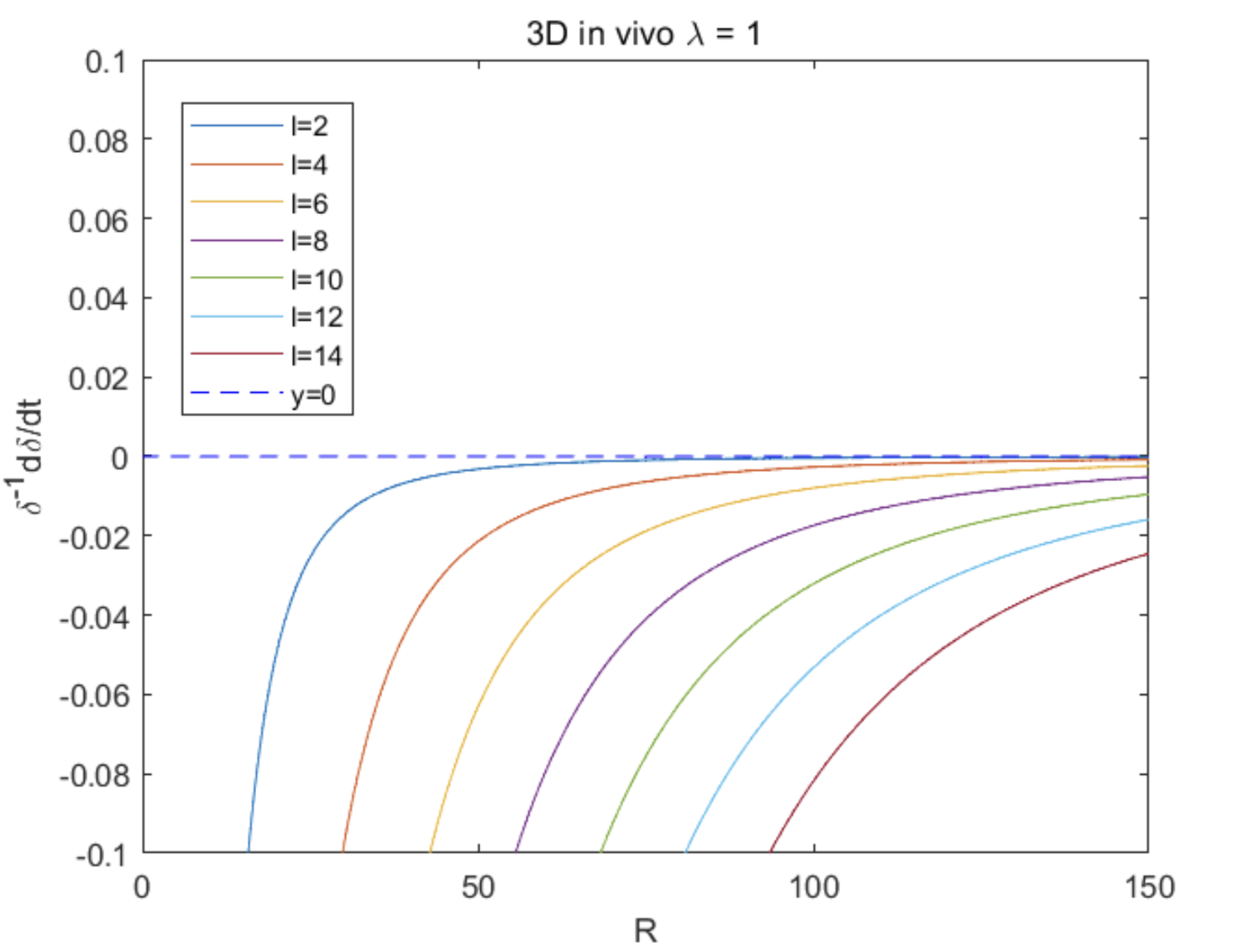}
        \caption{3D in vivo $\lambda=1,\ell=2,4,6,8,10,12,14, R\in[0,150]$} \label{vivo_lambda1_thres}
   \end{minipage}\hfill
   \begin{minipage}{0.53\textwidth}
        \centering
        \includegraphics[width=1.0\linewidth]{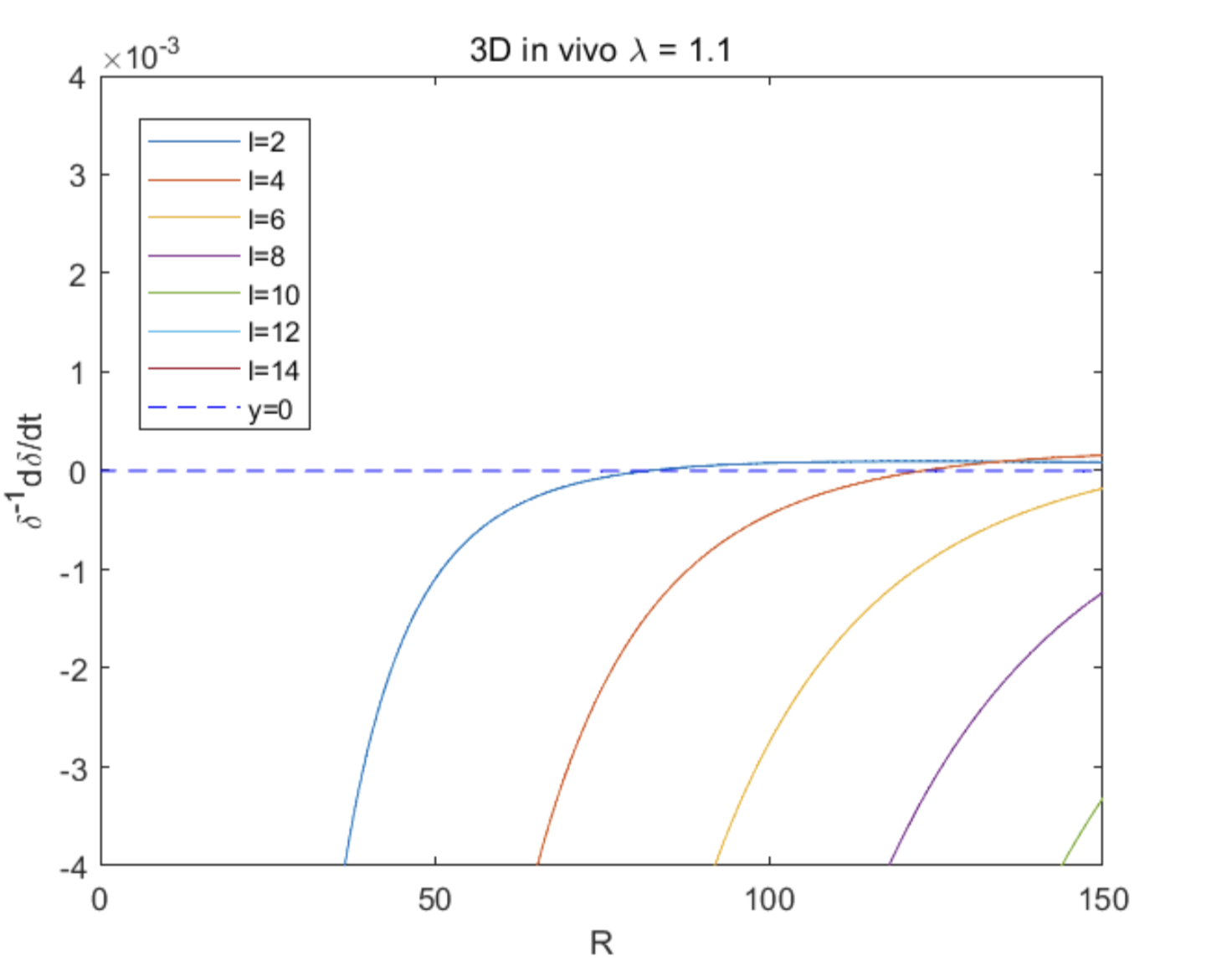}
        \caption{3D in vivo $\lambda=1.1,\ell=2,4,6,8,10,12,14, R\in[0,150]$}\label{vivo_lambda1.1_thres}
   \end{minipage}
   \end{figure}

\textbf{Theorem 5.1.} In 3D, given growing rate $G_0>0$, background concentration $c_B>0$, nutrient consumption rate $\lambda>0$, and perturbation degree $\ell$ and order $m$ of the spherical harmonic function $Y_{\ell}^{m}(\theta,\phi)$ as in (\ref{Ylm}). If the radius of the tumor is around $R$, where the corresponding unperturbed tumor has radius $R$, the evolution function $\delta^{-1}\frac{d\delta}{dt}$ for the in vitro model is given by Equation (\ref{3d vitro evol}):
\begin{equation} \nonumber
    \begin{aligned}      
        \delta^{-1}\frac{d \delta}{dt} &=  \frac{ R c_B G_0  I_{\ell+\frac{1}{2}}(\sqrt{\lambda}R) R^{-1-\ell}\sinh(\sqrt{\lambda}R) R^{\ell+\frac{5}{2}}\ell \sqrt{\lambda}}{\lambda^{3/2}R^{7/2}\sinh(\sqrt{\lambda}R) I_{\ell+\frac{1}{2}}(\sqrt{\lambda}R) }\\
        & +\frac{R c_B G_0  I_{\ell+\frac{1}{2}}(\sqrt{\lambda}R) (2\sqrt{\lambda}\sqrt{R}+\lambda^{3/2} R^{5/2})\sinh(\sqrt{\lambda}R)
        }{\lambda^{3/2}R^{7/2}\sinh(\sqrt{\lambda}R) I_{\ell+\frac{1}{2}}(\sqrt{\lambda}R)}  \\
        & -\frac{2R c_B G_0  I_{\ell+\frac{1}{2}}(\sqrt{\lambda}R) \cosh(\sqrt{\lambda}R)\lambda R^{3/2}}{\lambda^{3/2}R^{7/2}\sinh(\sqrt{\lambda}R) I_{\ell+\frac{1}{2}}(\sqrt{\lambda}R)}  \\
        &  -\frac{ (-\lambda(\ell+1)I_{\ell+\frac{1}{2}}(\sqrt{\lambda}R)+I_{\ell-\frac{1}{2}}(\sqrt{\lambda}R) \lambda^{3/2}R) R^{5/2} \cosh(\sqrt{\lambda}R) 
        } {2\lambda^{3/2}R^{7/2}\sinh(\sqrt{\lambda}R) I_{\ell+\frac{1}{2}}(\sqrt{\lambda}R) } +O(\delta) \\ 
        & \defeq M_1(\lambda,\ell,R)+O(\delta)
    \end{aligned}
\end{equation}
and for the in vivo model, it is given by (\ref{3d-vivo-delta}):
 \begin{equation} 
    \delta^{-1}\frac{d\delta}{dt} \defeq \frac{T_1}{T_2}+O(\delta) \defeq M_2(\lambda,\ell,R)+O(\delta)
\end{equation}
where,
\begin{equation} \label{T1}
    \begin{aligned}
        T_1& = - G_0 c_B( ((-2R( (3R+2)\lambda^{3/2}+\sqrt{\lambda}(R-2) )(\ell+1/2) K_{\ell+\frac{1}{2}}(R) \\
        &-2(R(R-\ell/2-1) \lambda^{3/2}-\sqrt{\lambda}(\ell+2)(R-2)/2)(R+1)K_{\ell-\frac{1}{2}}(R))\cdot  \\
        &\cosh{(\sqrt{\lambda}R)}^3-2\sinh{(\sqrt{\lambda}R)} ( R(\ell+1/2)( (\lambda^2+3\lambda) R+\lambda^2-1 ) K_{\ell+\frac{1}{2}}(R) \\
        &+((\lambda^2+\lambda)R^2-2\lambda(\ell+2)R+(\lambda+1)(\ell+2))(R+1)K_{\ell-\frac{1}{2}}(R)/2)\cdot\\ 
        & \cosh{(\sqrt{\lambda}R)}^2+(2R((2R+2)\lambda^{3/2}+\sqrt{\lambda}(R-2))(\ell+1/2)K_{\ell+\frac{1}{2}}(R) \\
        &+2(R^2\lambda^{3/2}-\sqrt{\lambda}(\ell+2)(R-2)/2)(R+1)K_{\ell-\frac{1}{2}}(R))\cosh{(\sqrt{\lambda}R)} \\
        & +2\sinh{(\sqrt{\lambda}R)}(R(\ell+1/2)(\lambda R+\lambda-1)K_{\ell+\frac{1}{2}}(R)\\
    \end{aligned}
\end{equation}   
\begin{equation} \nonumber
    \begin{aligned}
        &+K_{\ell-\frac{1}{2}}(R)(R+1)(R^2\lambda+\ell+2)/2))I_{\ell+\frac{1}{2}}(\sqrt{\lambda}R)\\
        &+I_{\ell-\frac{1}{2}}(\sqrt{\lambda}R)(R^2+(\ell+2)R+\ell+2)(\lambda(-2+(\lambda+1)R)\cosh{(\sqrt{\lambda}R)}^3\\
        &+2\sinh{(\sqrt{\lambda}R)}((R-1/2)\lambda^{3/2}-\sqrt{\lambda}/2)\cosh{(\sqrt{\lambda}R)}^2\\
        &-\lambda(R-2)\cosh{(\sqrt{\lambda}R)}+\sinh{(\sqrt{\lambda}R)}\sqrt{\lambda})K_{\ell+\frac{1}{2}}(R) \\
    \end{aligned}
\end{equation}
and,
\begin{equation} \label{T2}
    \begin{aligned}  
        T_2&= \sqrt{\lambda}R^3 K_{\ell+\frac{1}{2}}(R) (\cosh{(\sqrt{\lambda}R)}\lambda+ \sinh{(\sqrt{\lambda}R)}\sqrt{\lambda} )I_{\ell-\frac{1}{2}}(\sqrt{\lambda}R)\\
        &+I_{\ell+\frac{1}{2}}(\sqrt{\lambda}R) K_{\ell-\frac{1}{2}}(R)(\sqrt{\lambda} \cosh{(\sqrt{\lambda}R)} + \sinh{(\sqrt{\lambda}R)}))\cdot\\
        & (\cosh{(\sqrt{\lambda}R)}\lambda+ \sinh{(\sqrt{\lambda}R)}\sqrt{\lambda} )   (\sqrt{\lambda} \cosh{(\sqrt{\lambda}R)} + \sinh{(\sqrt{\lambda}R)})\\
    \end{aligned}
\end{equation}
Note that the scaling parameter $G_0, c_B>0$ do not influence the quantitative behavior of $\delta^{-1}\frac{d\delta}{dt}$. Also, in our analysis, the wavenumber $m$ is cancelled out in (\ref{Laplace c perturbed}). This is reasonable, considering that although $m$ primarily specifies orientations or symmetries, the fundamental characteristics of the solutions are predominantly determined by the parameter $\ell$. This is due to the dependency of $m$ on $\ell$, as indicated by the relationship $-\ell \leq m \leq \ell$. Essentially, it is wavenumber $\ell$ that plays the central role in defining the nature of these solutions.
For the in vitro model, we will show that $M_1(\lambda, \ell,R)$ is always negative, where they are shown in Figures \ref{3d vitro-1} and \ref{3d vitro-2}. For the in vivo model, fix the value of $G_0,c_B>0$, $M_2(\lambda, \ell,R)$ is plotted in Figures \ref{12}-\ref{15}, for different choice of $\lambda$ and perturbation parameter $\ell$. Based on the expression of  $\delta^{-1}\frac{d\delta}{dt}$ for the two nutrient models and the Figures, we establish the following remarks.\\ \\
\textbf{Remark 5.2.} $M_1(\lambda,1,R)=M_2(\lambda,1,R)=0$ for any $\lambda,R>0$. Since the model 1 perturbation corresponds to a trivial translation instead of the change of boundary geometry.\\  \\
\textbf{Remark 5.3.} 
When $0<\lambda<\lambda^{*}(\ell)$, fix any degree $\ell \leq 2$ and order $-\ell \leq m \leq \ell$ of the spherical harmonic function $Y_{\ell}^{m}(\theta,\phi)$ in the perturbation, $M_1(\lambda,1,R)$ is always negative and monotone increases in $R$ (see Figures \ref{12}-\ref{15}). From a physical perspective, when the nutrient consumption rate $\lambda$ is relatively low, the perturbation amplitude will continue to decrease to zero, where $\lambda^{*}(\ell)$ is the threshold of $\lambda$, and may be dependent on the wavenumber $\ell$. This means that the tumor in this range will always evolve from an irregular shape to a regular shape like sphere (with a larger size). \\ \\
\textbf{Remark 5.4.} 
For the regime $\lambda \geq \lambda^{*}(\ell)$, we have\\
(1) There exists a threshold $R^{*}(\ell)$ such that $F_2(\lambda,1,R)<0$ for $0<R<R^*(\ell)$, and $F_2(\lambda,1,R)>0$ for $R>R^*(\ell)$ (see Figures \ref{12}, \ref{13}). This means considering  degree $\ell$ and order $m$, and assume the nutrient consumption rate is significant, the perturbation amplitude will degenerate for tumor with small radius, and increase when the tumor boundary exceed the threshold radius ($R^*$), namely, the tumor will evolve from an irregular shape to a circular shape before the threshold radius $R^*$. Then it will continue to evolve after $R^*$, but from an spherical shape to an irregular 3D shape. \\ 
(2) Fix a proper value of $R_0$, there exists $\ell_0$ such that $M_2(\lambda, \ell,R)>0$ for $\ell<\ell_0$ and $M_2(\lambda, \ell,R)<0$ for $\ell>\ell_0$ (see Figures \ref{12}, \ref{13}), which implies that when the tumor size is around $R_0$, the perturbation of lower degree $\ell$ is easier to become unstable. However, the value of order $m$ is cancelled in (\ref{Laplace c perturbed}), and thus the threshold of tumor radius is independent on the order $m$. \\ 
(3) With the size of the tumor expanding, $R(t)$ exceeds more threshold $R^*(\ell)$, therefore the corresponding degree $\ell$ perturbation become unstable successively. \\ 
We hypothesize that there exists the threshold $\lambda^{*}$ between $1.0$ and $1.1$, as shown in Figures \ref{vivo_lambda1_thres} and \ref{vivo_lambda1.1_thres}, since the evolution function remains negative when $\lambda=1.0$ as in Figure \ref{vivo_lambda1_thres} and start going above $y=0$ when $\lambda=1.1$ as in Figure \ref{vivo_lambda1.1_thres}.

\textbf{Corollary 5.5.} Fix $G_0>0$ and $C_B>0$. For any $\lambda>0$ and $\ell>0$, $M_1(\lambda,\ell,R)$ Therefore, the perturbation amplitude always decay for the in vitro model:
\begin{equation}\label{M1}
\begin{aligned}
M_1(\lambda,\ell,R) &\sim \frac{-10 (\ell-1)}{3}+ O(R) \ \ \mbox{as} \ \ R\sim 0 \\
 M_1(\lambda,\ell,R) &\sim \frac{20\sqrt{\lambda} R}{(\ell+1)\left(-\frac{1}{4}\ell^2-\frac{1}{4} \ell +\frac{1}{2}\right)\ell }+O(1) \ \ \mbox{as} \ \ R\sim \infty 
 \end{aligned}
 \end{equation}
Therefore, $M_1(\lambda,\ell,R)<0$ as $R\to 0$ and as $R\to\infty$, thus $M_1(\lambda,\ell,R)<0$ always holds, since $\ell \geq 1$ is always true. 
Furthermore, for the in vivo model, we can infer from the figures that the function $M_2(\lambda,\ell,R)$ is always negative as $R\sim 0 $, and is positive as $R\sim \infty $. The existence of the threshold can be shown by the intermediate value theorem, where the function $M_2(\lambda,\ell,R)$ must intersect the horizontal axis, changing sign from negative to positive. This implies that the tumor boundary will change from stable to unstable. Note that the threshold $\lambda^{*}(\ell)$ may also depend on the degree $\ell$ as the perturbation parameter of spherical harmonics $Y_{\ell}^m$. These can be seem more clearly in Figure \ref{12}-\ref{15}.
In order to determine the critical threshold value $\lambda^{*}(\ell)$, we expect to substitute the asymptotic expansions of modified Bessel functions and hyperbolic functions into the evolution function as in (\ref{3d-vivo-delta}. The focus will be on the coefficient of the highest order term, as it dominates the behavior of $M_2(\lambda,\ell,R)$. By setting this leading coefficient to zero, we can ascertain the value of $\lambda^{*}$ that induces a sign transition in $M_2(\lambda,\ell,R)$. This transition, from negative to positive in the evolution function, marks the critical point where the stability of the tumor changes.

\subsection{2D and 3D results comparison}
In our comparative analysis of 2D and 3D tumor growth models, as initially discussed in \cite{feng2023tumor} and summarized in Section 2, we gained further insights into tumor boundary stability in different dimensions. In summary, the 2D and 3D models are essentially the same, and they both exhibit similar divergences between the in vitro and in vivo models. Particularly, for both the in vitro 2D and 3D models, the evolution function consistently remains negative, indicating stable tumor boundaries. 

For both the 2D and 3D in vitro models, the evolution function remains negative (indicating stability) when the consumption rate $\lambda$ is below the threshold $\lambda^{*}$, where $\lambda^{*}=1$ for the 2D model. Further investigation is required for the 3D model. When the rate exceeds this threshold, symbolized by $\lambda^{*}$, the function transitions from negative to positive, suggesting a shift from stable to unstable tumor boundaries. The threshold of evolution function $\lambda^{*}$ in 3D, and its possible dependence on the wavenumber $\ell$ will be further studied and compared with 2D scenario, beyong Figures (\ref{vivo_lambda1_thres}) and (\ref{vivo_lambda1.1_thres}).

The stability of tumor boundaries is thus determined by the sign of the evolution function: a negative evolution function indicates stable boundaries, whereas a positive evolution function signals instability. This distinction is crucial for understanding the dynamic behavior of tumor growth in different dimensional models, and the difference in 2D and 3D in vivo models mark the significance of our work.

\section{Conclusion}\label{sec13}
In conclusion, this study leverages numerical methods to validate and extend the analytical results presented in \cite{feng2023tumor} concerning 2D tumor PDE models and boundary instability. By reconstructing and analyzing main findings in \cite{feng2023tumor}, we successfully expanded these concepts into a 3D context. This extension is particularly significant as it offers a more applicable and insightful perspective on tumor behavior, especially regarding its instability and potential malignancy.
A key aspect of our approach is the combination of numerical simulations with analytical methods. This dual approach not only enriches our understanding of the tumor models but also serves as a cross-validation mechanism, ensuring the robustness of our findings.

During our research, we encountered two notable challenges in the numerical simulations, both of which provided significant insights. The first challenge arose when comparing the numerical tumor radius (\(R\)) over time against analytical speed (\(\sigma\)) over time in \cite{feng2023tumor}. The direct time derivative of numerical \(R\) data amplified errors, leading to impractical comparison results. We addressed this by integrating the analytical \(\sigma-t\) data for direct comparison with the numerical \(R(t)\) data, thereby minimizing error margins and aligning closely with analytical predictions in \cite{feng2023tumor}.
The second challenge involved the instability encountered in the ADI method when separately solving the reaction-diffusion equations for tumor density \(u\) and concentration \(c\). Although the ADI scheme's implementation was straightforward, its weakly decoupled system lacked the stability required for varied parameter inputs. To overcome this, we explored multiple splitting approaches within the ADI method, ultimately adopting the most stable splitting scheme.

In our final analysis, we examined the 3D in vivo evolution function by plotting it against various values of \(\ell\) and \(\lambda\). This led us to hypothesize the existence of a threshold, \(\lambda^{*}(\ell)\), which may or may not be influenced by the order \(\ell\). We observed that the evolution function remains negative for $\lambda<\lambda^{*}(\ell)$, and transits from negative to positive as $\lambda$ exceeds the threshold $\lambda^{*}(\ell)$. We propose that \(\lambda^{*}(\ell)\) could be determined using a methodology akin to the asymptotic analysis presented in Section 4.4, Step 5.
However, due to the limitations of numerical tools like Maple or Mathematica, currently it is hard to find the explicit formation of $\lambda^*(\ell)$. We save such formation for the future studies. Moreover, our intention encompasses extending beyond the current scope of linear stability analysis, which primarily focuses on the perturbation $\cos(m\phi)$ in 2D and $Y_{\ell}^m(\theta, \phi)$ in 3D. Our planned trajectory involves incorporating nonlinear terms, aligning with our prospective research objectives.


\backmatter

\bmhead{Supplementary information} Relevant code are available at:
\url{https://github.com/JiaqiZhang988/Tumor-Growth-PDE-Model-Relevant-Code}

\bmhead{Acknowledgments}
This project is partially supported by National Key R\&D Program of China, Project Number 2021YFA1001200. J. Zhang is partially supported by Summer Research Scholar program from Duke Kunshan University. X. Xu is partially supported by Natinal Science Foundation of China Youth program, Project number 12101278, and Kunshan Shuangchuang Talent Program, Project Number kssc202102066. The authors would also like to thank the helpful discussion with Yu Feng, Dang Xing Chen and Lin Jiu.

\section*{Declarations}
\begin{itemize}
\item Conflict of interest \\
On behalf of all authors, the corresponding author states that there is no conflict of interest.
\item Funding\\
This project is partially supported by National Key R\&D Program of China, Project Number 2021YFA1001200. 
\item Ethics approval \\
Not applicable.
\item Consent to participate\\
Not applicable.
\item Consent for publication\\
Not applicable.
\item Availability of data and materials\\
Data and supporting materials for this work are available from the corresponding author upon reasonable request.
\item Code availability \\
Code of this work are available at \url{https://github.com/JiaqiZhang988/Tumor-Growth-PDE-Model-Relevant-Code}.
\item Authors' contributions\\
This paper was completed by J. Zhang under the guidance and supervision of Prof. T. Witelski, Prof. J. Liu and Prof. X. Xu. All authors contributed to refining the methodologies and to the writing and reviewing of the manuscript. Each author has read and approved the final version of the manuscript for publication.

\end{itemize}



\begin{appendices}
\section{Asymptotic expansions of modified Bessel functions and hyperbolic trigonometric functions} \label{secAAA}
When $r \rightarrow 0$ \cite{NIST_DLMF_Bessel}:
\begin{equation} \label{asym-IK}
    \begin{aligned}
        \mbox{when} \ \ r \rightarrow 0, \\
        I_{\nu}(r) \sim \sum_{k=0}^{\infty}\frac{1}{\Gamma(k+1)\Gamma(k+\nu+1)}\left( \frac{r}{2} \right)^{2k+\nu}, \ \ \mbox{for} \ \ \nu \in \mathbb{R},\\
         K_{\nu}(r) \sim \frac{\pi}{2}\frac{I_{-\nu}(r)-I_{\nu}(r)}{\sin(\nu\pi)} , \ \ \mbox{for} \ \ \nu \in \mathbb{R}.
    \end{aligned}
\end{equation}
and 
\begin{equation} \label{asym-coshsinh}
    \begin{aligned}
         \mbox{when} \ \ r \rightarrow 0, \\
         \cosh(r) \sim 1+\frac{r^2}{2!}+\frac{r^4}{4!}+\frac{r^6}{6!}+\frac{r^8}{8!}+\frac{r^{10}}{10!}+\cdots \\
         \sinh(r) \sim r+\frac{r^3}{3!}+\frac{r^5}{5!}+\frac{r^7}{7!}+\frac{r^9}{9!}+\frac{r^{11}}{11!}+\cdots
    \end{aligned}  
\end{equation}
When $(r \rightarrow \infty)$:
\begin{equation}
    \begin{aligned}
      \mbox{when} \ \ r \rightarrow \infty, \\
        I_{\nu}(r) \sim \frac{e^r}{(2\pi r)^\frac{1}{2}}\sum_{k=0}^{\infty}(-1)^k \frac{a_k(\nu)}{r^k} \\
        K_{\nu}(r) \sim \left( \frac{\pi}{2r}  \right)^{\frac{1}{2}}e^{-r} \sum_{k=0}^{\infty}\frac{a_k(\nu)}{r^k} \\
        \mbox{where}  \ \ \ a_k(\nu)=\frac{(4\nu^2-1^2)(4\nu^2-3^2)\cdots(4\nu^2-(2k-1)^2)}{k!8^{k}}\\   
    \end{aligned}
\end{equation}
and 
\begin{equation}
    \begin{aligned}
     \mbox{when} \ \ r \rightarrow \infty, \\
        \cosh(r) \sim \frac{e^r}{2}\\
        \sinh(r) \sim \frac{e^r}{2}
    \end{aligned}
\end{equation}

\section{Two-Dimensional (2D) Numerics}\label{secA1}
\subsection{Numerical methods introduction}
To implement our numerical scheme in 2D model, we describe the $\Psi$ function for the in vivo regime as: 
\begin{equation}
    \Psi(\rho,c) = \Biggl\{ 
    \begin{aligned}
        \lambda \rho c \ \ \ \ \  \ \ \ \ \ \mbox{inside} \ \ \mbox{tumor} \ (\rho>0) \\
         -(c_B-c) \ \ \ \ \ \ \mbox{outside} \ \ \mbox{tumor} \ (\rho=0)
\end{aligned}
\end{equation}
Combining these, we define a $H(\rho)$ function to represent the function $\Psi(\rho,c)$ in numerical simulations, where
\begin{equation}
    \Psi(\rho,c)=-(c_B-c)+(\lambda \rho c +(c_B-c)) H(\rho)
\end{equation}
and we define the function $H(u)$ using the hyperbolic tangent function
\begin{equation}
    \begin{aligned}
        H(\rho) = \Biggl\{
        \begin{aligned}
              0, \ \ \ \mbox{if} \ \ \rho\leq 0 \\
              \mbox{tanh}(ku), \ \ \mbox{if} \ \ \rho>0. \\
    \end{aligned}  \\
    H'(\rho) = \Biggl\{
    \begin{aligned}
         0, \ \ \ \mbox{if} \ \ \rho< 0 \\
        k \mbox{sech} ^2 (ku), \ \ \ \mbox{if} \ \ \rho>0.
    \end{aligned}
\end{aligned}
\end{equation}
where $k$ is a constant, and the bigger the $k$ is, the narrower the smooth part of the step function is. Including the hyperbolic function with the constant $k$ gives us an approximation of the smoothed step function and helps us avoid the sharp change of the step function when using the indicator function $\mathds{1}$, thus avoiding problems brought by the steep step function in the numerical process. \\
\subsubsection{2D ADI model.} The general PDEs for the 2D reaction diffusion system based on Equation(\ref{govern}) and (\ref{govern c}), after rewriting the PME part (with power $u^k$), we have
\begin{equation} 
    \begin{aligned}
        {\rho}_t = \frac{k}{k+1}( ({\rho}^{k+1})_{xx} +  ({\rho}^{k+1})_{yy}) + G_0 c \rho \\
        c_t = \frac{1}{\tau} (c_{xx} + c_{yy} - \Psi(\rho,c))
    \end{aligned}  
\end{equation}
Therefore, using finite difference and rewriting the discrete points above as ${\rho}_{i+1,j}$, below ${\rho}_{i-1,j}$, to the right ${\rho}_{i,j+1}$ and to the left $u_{i,j-1}$ of the center point ${\rho}_{i,j}$, where $i$ is the index in the $x$ direction, and $j$ is the index in the $y$ direction:
\begin{equation}
    \begin{aligned}
        {\rho}_{i,j} = {\rho}_c \\ 
        {\rho}_{i+1,j} = {\rho}_u, \ \  {\rho}_{i-1,j}={\rho}_d, \\
        {\rho}_{i,j+1} = {\rho}_p, \ \ {\rho}_{i,j-1}={\rho}_m
    \end{aligned}
\end{equation}
The same applies to $c$, where we have 
\begin{equation}
    \begin{aligned}
        c_{i,j} = c_c \\ 
        c_{i+1,j} = c_u, \ \  c_{i-1,j}=c_d, \\
        c_{i,j+1} = c_p, \ \ c_{i,j-1}=c_m
    \end{aligned}
\end{equation}
we have the discretized version of this 2D reaction diffusion system  
\begin{equation} \label{2d ADI r-d discretized}
    \begin{aligned}
        \frac{{\rho}_c -{\rho}_{c,old}}{\Delta t} = &\frac{k}{k+1} \left(  \left( \frac{{\rho}_p^{k+1} - 2{\rho}_c^{k+1} + {\rho}_k^{k+1} }{\Delta x^2} \right) +  \left( \frac{{\rho}_u^{k+1} - 2{\rho}_c^{k+1} + {\rho}_d^{k+1} }{\Delta y^2} \right) \right)\\
         &+ G_0 c_i {\rho}_i  \\
        \frac{c_c-c_{c,old}}{\Delta t} =& \frac{1}{\tau} \left( \left( \frac{c_p-2c_c+c_m}{\Delta x^2}  \right) +  \left( \frac{c_u-2c_c+c_d}{\Delta y^2}  \right) \right) - \frac{1}{\tau}\Psi({\rho}_i,c_i)
    \end{aligned}
\end{equation}
where ${\rho}_c$ and ${\rho}_{c,old}$ represents the $u(x,y,t)$ value of the $i,j$th discretized grid point at the current time and the previous time, respectively. We then multiply both sides by $\Delta t$ and move all terms to the same side, and use them as the functions where we apply Newton's method finding its zeros. 

The initial condition (IC) for density $u$ takes the 2d axisymmetric form of the similarity solution of PME:
\begin{equation} \label{2d simi}
    {\rho}(t=0,r) = \left(  1-\frac{k}{4(k+1)^2} (r(t+1)^{\frac{-1}{2k+2}})^2 \right)^{\frac{1}{k}} (t+1)^{\frac{-1}{k+1}}
\end{equation}

The IC for concentration $c$ is the constant $c_B=0.1$ in this 2D case. The boundary conditions (BCs) we used are still Neumann BCs at both ends of the $x$ and $y$ directions, that is
\begin{equation}
    \begin{aligned}
        \frac{\partial {\rho}}{\partial r} \big|_{x=0} &=& 0,
        \frac{\partial c}{\partial r} \big|_{x=0} &=&0 ,
        \frac{\partial {\rho}}{\partial r} \big|_{x=R} &=&0 , 
        \frac{\partial c}{\partial r} \big|_{x=R} &=&0    \\ \nonumber
         \frac{\partial {\rho}}{\partial r} \big|_{y=0} &=&0, 
         \frac{\partial c}{\partial r} \big|_{y=0} &=&0,
         \frac{\partial {\rho}}{\partial r} \big|_{y=R} &=&0,
         \frac{\partial c}{\partial r} \big|_{y=R} &=&0  
    \end{aligned} 
\end{equation}

The parameters remain the same as the 1D scenario, namely, the growth constant  $G_0=1$, the consumption rate $\lambda=0.5$ (involved
in $\Psi$ term), and the characteristic time scale of nutrients evolution $\tau = 0.1$, and the constant for
the approximated smoothed step function $k=5$. 
\subsubsection{2D Axisymmetric Model.} 
We also consider the 2D axisymmetric model for comparison. Since the Laplacian for 2D is 
\begin{equation} \label{2D Laplacian}
    \Delta = \frac{1}{r}\frac{\partial {\rho}}{\partial r} \left( r \frac{\partial({\rho}^{k+1})}{\partial r}  \right)
\end{equation}
thus the general PDEs for 2D axisymmetric model are
\begin{eqnarray}
    {\rho}_t &=& \frac{k}{k+1} \frac{1}{r}\frac{\partial {\rho}}{\partial r} \left( r \frac{\partial({\rho}^{k+1})}{\partial r}  \right) + G_0c{\rho}  \\
    c_t &=& \frac{1}{\tau} \left(  \frac{1}{r}\frac{\partial c}{\partial r} \left( r \frac{\partial(c)}{\partial r}  \right) - \Psi({\rho},c)  \right)  \nonumber
\end{eqnarray}
Therefore, using finite difference method and denote $u$ at the center point as ${\rho}_i$ with discretization step $\Delta r$ and $\Delta t$, we have
\begin{equation} \label{2d r-d discretized}
    \begin{aligned}
         \frac{{\rho}_i - {\rho}_{i,old}}{\Delta t} =& \frac{k}{k+1}\frac{1}{r_i \Delta r}\left[  \frac{ (r_i\!+\!0.5 \Delta r) ({\rho}_{i+1}^{k+1}-{\rho}_i^{k+1}) }{\Delta r}  -  \frac{ (r_i-0.5 \Delta r) ({\rho}_{i}^{k+1}-{\rho}_{i-1}^{k+1}) }{\Delta r}   \right]  \\
         & + G_0 c_i {\rho}_i\frac{c_i-c_{i,old}}{\Delta t} \\
         =& \frac{1}{\tau} \frac{1}{r_i \Delta r} \left[  \frac{ (r_i+0.5 \Delta r) (c_{i+1}-c_i) }{\Delta r}  -  \frac{ (r_i-0.5 \Delta r) (c_{i}-c_{i-1}) }{\Delta r}   \right] -\frac{1}{\tau}\Psi({\rho}_i,c_i) 
    \end{aligned}
\end{equation}
where ${\rho}_i$ and ${\rho}_{i,old}$ represents the $u(x,t)$ value of the $i$th discretized grid point at the current time and the previous time, respectively. We then multiply both sides by $\Delta t$ and move all terms to the same side, and use them as the functions where we apply Newton's method finding its zeros. 

The IC for density $u$ and concentration $c$ are the same as the 2D ADI model, and the parameters are also the same. The BCs are the same, where we have Neumann boundary conditions on both ends:
\begin{eqnarray} \label{2d axis bcs}
    \frac{\partial {\rho}}{\partial r} \Big|_{r=0} &=& 0, \ \ \
    \frac{\partial c}{\partial r} \Big|_{r=0} =0 \\ \ \nonumber
    \frac{\partial {\rho}}{\partial r} \Big|_{r=R} &=& 0, \ \ \ \nonumber
    \frac{\partial c}{\partial r} \Big|_{r=R} =0   \nonumber
\end{eqnarray}

\subsection{Numerical methods implementation.}
\subsubsection{2D ADI for PME}
In this section, we follow the ADI scheme introduced in \cite{witelski2003adi} to solve the 2D PME and the reaction-diffusion equation systems in \cite{feng2023tumor}. Specifically, we use the approximate factorization and numerically solve the 2D nonlinear PDE problem by solving 1D nonlinear problem twice, namely, first in $x$, and then in $y$ direction.

We consider the porous medium equation (PME) in 2D:
\begin{equation}
    {\rho}_t = \nabla^2 ({\rho}^{k+1})
\end{equation}

Using backward Euler implicit method in finite difference method, we can discretize the equation
\begin{equation}
    \frac{{\rho}-{\rho}^{old}}{\Delta t} = \nabla^2 ({\rho}^{k+1}) + O(\Delta t)
\end{equation}

To use the Newton's iteration method to solve this implicit scheme, we write
\begin{equation} \label{F_adi}
    F({\rho}) = {\rho}-\Delta t \nabla^2 ({\rho}^{k+1}) - {\rho}^{old} + O(\Delta t^2)
\end{equation}

We set ${\rho}_0={\rho}^{old}$ and write $F(u)=0$ as
\begin{equation}
    \begin{aligned}
        J_k \Delta {\rho}_k = -F_k \\
        {\rho}_{k+1} = {\rho}_k + \Delta {\rho}_k
    \end{aligned}
\end{equation}

Here, $J_{k_1} \Delta {\rho}_{k_1} = \frac{\delta F}{\delta {\rho}} \Delta u $ is the functional derivative of $F$ (the Jacobian), which we approximate by $F({\rho}+\Delta {\rho}) - F(u)$, plug in ${\rho}+\Delta {\rho}$ and ${\rho}$ into Eq.~\ref{F_adi} and only keep the linear terms in $\Delta {\rho}$ and neglect the higher order terms, we can arrive at:
\begin{equation}
    \begin{aligned}
        J_{k_1} \Delta {\rho}_{k_1} &= F({\rho}+\Delta {\rho}) - F({\rho})\\
        &= \Delta {\rho}_{k_1} - \Delta t \nabla^2((k+1){\rho}_{k_1}^k \Delta {\rho}_{k_1}) \\
        &= \Delta {\rho} - \Delta t (k+1)({\rho}_{k_1}^k \Delta {\rho})_{xx} - \Delta t(k+1)({\rho}_{k_1}^k \Delta {\rho})_{yy}
    \end{aligned}
\end{equation}

The matrix operators following this are: ($I$: identity matrix)
\begin{equation}
    (I - \Delta t D_{xx} - \Delta t D_{yy}) \Delta {\rho}
\end{equation}
where we denote $D_{xx} = (k+1)({\rho}_{k_1}^k \Delta {\rho})_{xx} $, and $D_yy = (k+1)({\rho}_{k_1}^k \Delta {\rho})_{yy}$. 

Therefore, we can approximate the 2D Jacobian matrix $J$ via
\begin{equation}
    \begin{aligned}
        (I- \Delta t D_{xx})(I- \Delta t D_{yy}) \Delta u \\
        (I-\Delta t D_{xx} - \Delta t D_{yy} + O(\Delta t^2)) \Delta u    
    \end{aligned}
\end{equation}
where can ignore the order $O(\Delta t^2)$ term, which means that we can use $(I-\Delta t D_{xx} - \Delta t D_{yy})$ as an approximation of $J_{k_1} \Delta u_{k_1}$, namely, $(I-\Delta t D_{xx} - \Delta t D_{yy}) \approx Jk \Delta u_k$. 

Therefore, we have
\begin{equation}
    (I- \Delta t D_{xx})[(I- \Delta t D_{yy}) \Delta {\rho}] = -F_{k_1} 
\end{equation}

let $v = (I- \Delta t D_{yy}) \Delta {\rho}$, we have
\begin{equation}
    \Biggl\{ 
    \begin{aligned}
        (I- \Delta t D_{xx}) v = -F_{k_1} \\
        (I- \Delta t D_{yy}) \Delta {\rho} = v
    \end{aligned}
\end{equation}

We first solve the $D_{xx}$ equation to get $v$, and then substitute $v$ into the $D_{yy}$ equation to solve for $\Delta {\rho}$, which will be used to update the initial guess in the Newton iteration, and obtain improved approximation through each update until the solution converges.

Specifically, we first solve the $D_{xx}$ equation for $v$ (substitute in expression for $D_{xx}$:
\begin{equation}
    \begin{aligned}
        (I- \Delta t D_{xx}) v = -F_{k_1} \\
        \Rightarrow v - \Delta t D_{xx} v = -F_{k_1} \\
        \Rightarrow v - \Delta t (k+1)({\rho}^k v)_{xx} = -F_{k_1} \\
        \Rightarrow v_i - \Delta t (k+1) \left( \frac{{\rho}_{i+1}^k v_{i+1} - 2{\rho}_i^k v_i + {\rho}_{i-1}^k v_{i-1} }{\Delta x^2} \right) = -F_{k_1}       
    \end{aligned}
\end{equation}

This can be expressed in matrix form:
\begin{equation}
   \left[ I - (k+1)\frac{\Delta t}{\Delta x^2}
    \begin{pmatrix}
        -2{\rho}_i^k & {\rho}_{i-1}^k & & & & \\
        {\rho}_{i+1}^k & 2{\rho}_i^k & {\rho}_{i-1}^k & &  &\\
        &\ddots & \ddots & \ddots &\\
        & & {\rho}_{i+1}^k & 2{\rho}_i^k & {\rho}_{i-1}^k &\\
        & & & -2{\rho}_i^k & {\rho}_{i+1}^k & \\ 
    \end{pmatrix} \right] 
    \begin{pmatrix}
        v_1\\
        v_2\\
        \vdots\\
        v_M
    \end{pmatrix}
    = 
    \begin{pmatrix}
    F_1\\
    F_2\\
    \vdots\\
    F_M
    \end{pmatrix}
\end{equation}

Using backslash in MATLAB, we can get vector $v$. Note that for this first step, we solve in the direction of $x$, we take each row of matrix $F$, and solve the corresponding vector $v$ for each row of matrix $F$. All vector $v$ together forms a matrix $v$. Then we take use similar method to calculate $\Delta u$ in the $y$ direction. Note that we take each column of matrix $v$, and use the matrix form $(I-\Delta t D_{yy})\Delta {\rho} = v$ to calculate corresponding $\Delta {\rho}$ for each column of $v$. Finally, we add $\Delta {\rho}$ on ${\rho}_0$ to update initial guess and derive the numerical result of solution ${\rho}$ of PME using 2D ADI scheme.

In summary, the key idea in ADI scheme is to approximate the big Jacobian matrix $(J)$ for the 2D PME by factorization in $x$ and $y$ direction, namely, we factorize the Jacobian matrix $J=J_x J_y$, where $J_x = I-\Delta t D_{xx}$, and $J_y = I-\Delta t D_{yy}$.  
\subsubsection{2D axisymmetric for PME} 
As introduced in the general 2D axisymmetric model, we have the finite difference discretization for the PME part, which is
\begin{equation}
    \frac{k}{k+1}\frac{1}{r_i \Delta r} \left[  \frac{ (r_i+0.5 \Delta r) ({\rho}_{i+1}^{k+1}-{\rho}_i^{k+1}) }{\Delta r}  -  \frac{ (r_i-0.5 \Delta r) ({\rho}_{i}^{k+1}-{\rho}_{i-1}^{k+1}) }{\Delta r}   \right]
\end{equation}
where the IC is the similarity solution for 2D axisymmetric version of PME as in (\ref{2d simi}), and the boundary conditions (BCs) are the Neumann BCs at both ends, as given in (\ref{2d axis bcs}). 

\begin{figure}[!htb]
   \begin{minipage}{0.51\textwidth}
        \centering
        \includegraphics[width=1.0\linewidth]{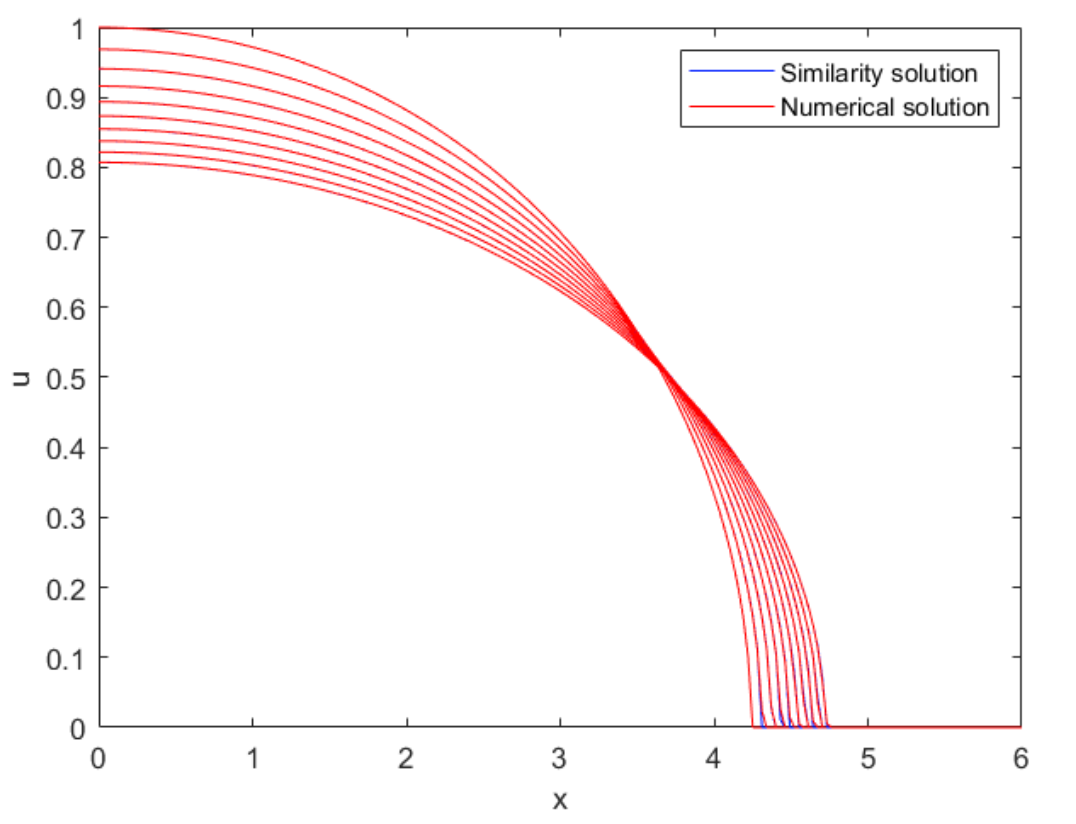}
        \caption{2D axisymmetric PME-Comparison of numerical and analytical solution}\label{Fig:2d axis-1}
   \end{minipage}\hfill
   \begin{minipage}{0.507\textwidth}
        \centering
        \includegraphics[width=1.0\linewidth]{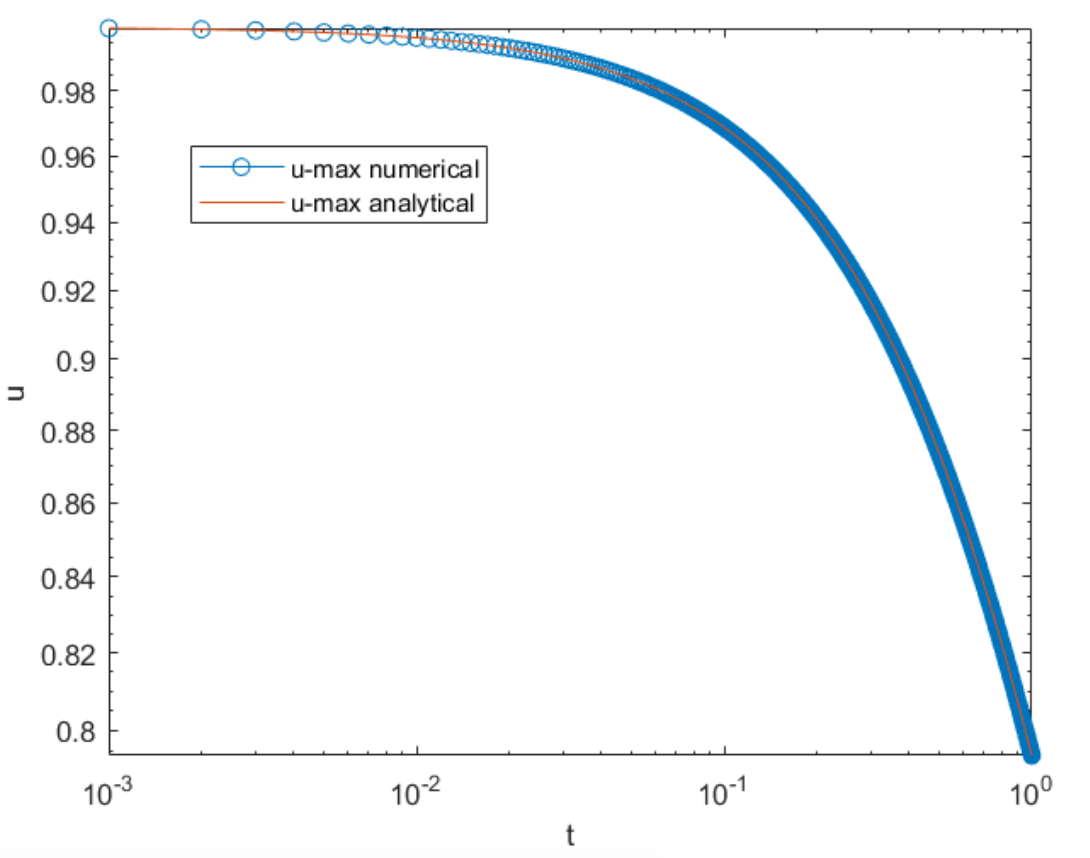}
        \caption{validation of 2D axisymmetric PME}\label{Fig:2d axis-2}
   \end{minipage}
\end{figure}

The parameters are: $L=20$, $T=50$, $M=100$, $N=1000$, $k=2$, where $L$ is the total length, $T$ is the final time, $M$ is the number of grid points in spatial direction, $N$ is the number of grid points in time discretization, and $k$ is the power of the porous medium equation. We use Newton's method to obtain the approximated solution that we believe converge fast and well to the true solution. The check of this is achieved by drawing the plot of the similarity and numerical solutions in the same figure for comparison (Figure.\ref{Fig:2d axis-1}). We also draw the log-log plot of the maximum of the approximated solution and the log-log plot of $(1+t)^{-1/4}$ in the same figure for comparison (Figure.\ref{Fig:2d axis-2}). 
\subsubsection{2D ADI and 2D axisymmetric comparison.}
In this section, we consider the porous medium equation (PME), and match up the two-dimensional (2D) ADI method with the 2D axisymmetric case, thus proving that the ADI method is valid, and we can have irregular shape or add perturbation ($\theta$ dependence) to be the initial condition and investigate the stability of the boundary after some time, namely, simulating the perturbation analysis and the instability tumor boundary in 2D and match up to \cite{feng2023tumor}. Specifically, we try to compare with the 2d axisymmetric similarity solution starting from an initial condition (IC) in the center of a square ADI domain. Then we start with a different IC, like a small ellipse, and we saw that it becomes more and more circular.  
\subsubsection{Method for matching up.} Firstly, we want to find the center of the initial shape, so that we can allow MATLAB easily find the center of any irregular shape. This can be achieved by utilizing the mass conservation of PME with Neumann-Neumann boundary conditions. 
Specifically, for ${\rho}_t=\nabla^2 ({\rho}^k)$ with ${\rho}_x=0$, we have
\begin{equation}
    \int {\rho} \, dA\ = \mbox{constant}, \ \ \int x{\rho} \, dA\ = \mbox{constant}, \ \ \int y{\rho} \, dA\ = \mbox{constant}
\end{equation}
and thus, the center ($\bar{x}, \bar{y}$) is
\begin{equation}
    \bar{x} = \frac{\int x{\rho} \, dA\ }{\int {\rho} \, dA\ }, \ \ 
    \bar{y} = \frac{\int y{\rho} \, dA\ }{\int {\rho} \, dA\ }  
\end{equation}
In this way, we can do this calculation and get the position of the x and y center point calculated from the last time, so that we can define the center of any strange irregularly shaped tumor. 

The top view of the result of ADI is the circular region, thus we can go from that center point and go to four radial directions of the circular region. Either one of these four directions can be used to check against the 2D axisymmetric case. Specifically, for the ADI we use one of the directions and find its corresponding ${\rho}(r)$ with $r$ going from the center to the end point, and draw the ${\rho}(r)$ versus $r$ plot. We then compare this to the ${\rho}(r)$ versus $r$ plot generated at the final time in 2D axisymmetric case. We can see that the 2D ADI plot and the 2D axisymmetric case match up to be the same with the same initial condition (IC). 

This means that the 2D ADI can work to simulate thew 2D axisymmetric tumor considered in \cite{feng2023tumor}. The advantage of this is that for axisymmetric case, we can not have any $\theta$ dependence and only depends on $r$, but for ADI, we can let the solution to evolve in $\theta$. This allows us to simulate the stability of the tumor boundary after some perturbation dependent on $\theta$. Specifically, if we start from IC $ {\rho}_0 = f(r)$, then for the tumor model, we will have some $ \hat{{\rho}_0} = f(r) + \eps g(\theta)$ (small perturbation depend on $\theta$), which lead to small change in the shape of the tumor. With the ADI code, we can observe the stability of the tumor, namely, does the tumor go back to being almost circular, or does the perturbation grow and the tumor shape changes from being nice to being more messy looking. Note that we also check that the ADI scheme by plotting the $U_{max}$ versus time ($t$) and the $t^{-1/(k+1)}$ versus time ($t$) in the same $\log\log$ plot and we can see that they fit together to be a straight line with negative slope and are the same to each other. This is because the maximum of ${\rho}_{max}(r)=t^{-1/(m+1)}$ and thus $\log ({\rho}_{max}(r)) = (-1/(k+1))\log t$. 

\begin{figure}
    \centering
    \includegraphics[width=8.2cm]{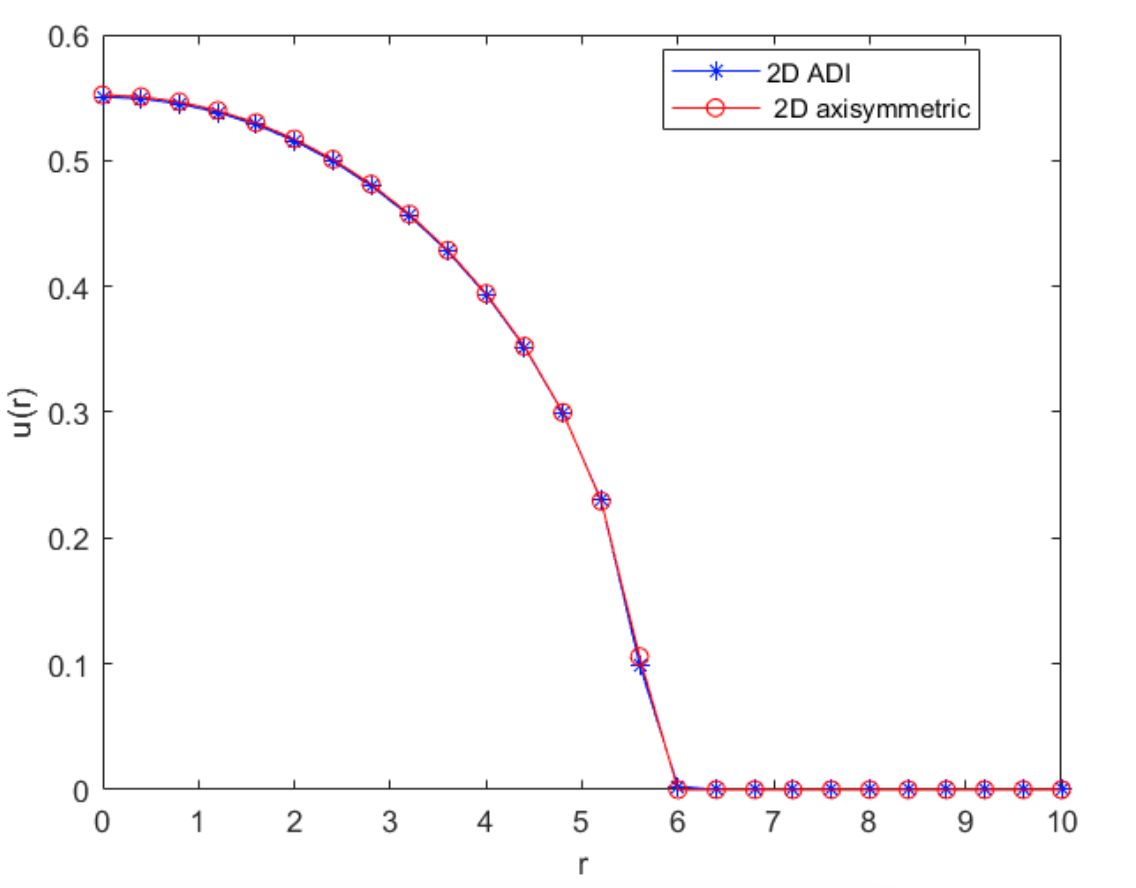}
    \caption{one plot: Comparison of 2D ADI and axisymmetric}
    \label{fig:adi and 2d axis compare}
\end{figure}

The comparison are shown in Figure\ref{fig:adi and 2d axis compare}, and we can see that the results from 2D ADI and 2D axisymmetric match up. 
\subsection{Validating Numerical simulations against Analytical Results and Predictions in \cite{feng2023tumor}}
\subsubsection{2D ADI reaction-diffusion systems numerical check for analyitcal results in \cite{feng2023tumor}:}
We implement the 2D ADI scheme on the reaction-diffusion systems to numerically check against the analytical results in \cite{feng2023tumor}. The discretization for the 2D axisymmetric reaction diffusion system is given in (\ref{2d ADI r-d discretized}). For the boundary conditions, we use the ghost point method to included the boundary conditions in the $F_u$ and $F_c$ equations, 
where,
\begin{equation}
    \begin{aligned}
        F_{{\rho}_i} = & {\rho}_c -{\rho}_{c,old} - \frac{k\Delta t }{k+1} \left(  \left( \frac{{\rho}_p^{k+1} - 2{\rho}_c^{k+1} + {\rho}_k^{k+1} }{\Delta x^2} \right) +  \left( \frac{{\rho}_u^{k+1} - 2{\rho}_c^{k+1} + {\rho}_d^{k+1} }{\Delta y^2} \right) \right)\\
        &- \Delta t G_0 c_i {\rho}_i  \\
        F_{c_i} = &c_c-c_{c,old} - \frac{\Delta t}{\tau} \left( \left( \frac{c_p-2c_c+c_m}{\Delta x^2}  \right) +  \left( \frac{c_u-2c_c+c_d}{\Delta y^2}  \right) \right) + \frac{\Delta t}{\tau}\Psi({\rho}_i,c_i)
   \end{aligned}
\end{equation}
Then we use $F_u, F_c$ to apply Newton's method as introduced before. \\ \\  
\textbf{$\delta^{-1} \frac{d \delta}{dt}$ check for the in vivo model.} \\
We also use the numerical simulation to check the perturbation analysis in \cite{feng2023tumor}. The analytical results for the tumor boundary instability after perturbation are: For the in vivo model, when the consumption rate $0<\lambda \leq 1$, if we fix any wave number $l \geq 2$ in the perturbation, then the instability will degenerate and the tumor boundary will become stable. Hence, the shape will change from irregular to circular with $0<\lambda \leq 1$. On the other hand, if the consumption rate $\lambda>1$, if we fix any wave number $l \geq 2$ in the perturbation, then there is a threshold radius $R^*$. This means that when the tumor radius $R \leq R^*$, then the tumor is stable, and after the radius exceed $R^*$, the tumor will eventually become unstable. Therefore, the shape will change from irregular to circular and then become irregular again with $\lambda>1$. The results showing the shape evolution of the tumor boundary is shown in Figures.\ref{3},\ref{4} in section 2.




\end{appendices}

\bibliography{bibliography} 

\end{document}